\documentclass[aos,MSNbibl,nameyear,seceqn,dvips]{arximspdf}
\usepackage{graphicx}

\doi{10.1214/12-AOS1075} %kopijuoti is PTS
\volume{41}
\issue{1}
\pubyear{2013}
\firstpage{91}
\lastpage{124}

\makeatletter
\newcommand{\rrvert}{\vert}
\newcommand{\llvert}{\vert}
\newcommand{\N}{\mathbb{N}}
\newcommand{\R}{\mathbb{R}}

\newcommand{\X}{\mathcal{X}}
\newcommand{\E}{\mathcal{E}}
\newcommand{\iidsim}{\stackrel{\mathrm{i.i.d.}}{\sim}}
\newcommand{\btheta}{\bolds{\theta}}
\newcommand{\bbeta}{\bolds{\beta}}
\newcommand{\balpha}{\bolds{\alpha}}
\newcommand{\bA}{\mathbf{A}}
\newcommand{\bS}{\mathbf{S}}
\newcommand{\bC}{\mathbf{C}}
\newcommand{\bx}{\mathbf{x}}
\newcommand{\by}{\mathbf{y}}
\newcommand{\bz}{\mathbf{z}}
\newcommand{\bN}{\mathbf{N}}
\newcommand{\cX}{\mathcal{X}}
\newcommand{\cY}{\mathcal{Y}}
\newcommand{\Gap}{\operatorname{Gap}}
\newcommand{\indic}{\mathbf{1}}
\newtheorem{prop}{Proposition}[section]
\newtheorem{theorem}{Theorem}[section]
\newtheorem{lemma}{Lemma}[section]
\newtheorem{cor}{Corollary}[section]
\newproclaim{assumption}{Assumption}[section]
\makeatother

\begin{document}
\begin{frontmatter}

\title{Convergence rate of Markov chain methods for genomic motif discovery\thanksref{T1}}
\runtitle{Convergence rate of Markov chains}

\thankstext{T1}{Supported in part by U.S. NSF awards CMMI-0926814 and DMS-12-09103 and by NSERC of Canada.}

\begin{aug}
\author[A]{\fnms{Dawn B.} \snm{Woodard}\corref{}\ead[label=e1]{woodard@cornell.edu}\ead[label=u1,url]{http://people.orie.cornell.edu/woodard}}
\and
\author[B]{\fnms{Jeffrey S.} \snm{Rosenthal}\ead[label=u2,url]{http://www.probability.ca/jeff}}
\runauthor{D. B. Woodard and J. S. Rosenthal}
\affiliation{Cornell University and University of Toronto}
\address[A]{School of Operations Research\\
\quad and Information Engineering\\
\quad and Department of Statistics\\
Cornell University\\
206 Rhodes Hall\\
Ithaca, New York 14853\\
USA\\
\printead{u1}} %adresu isvedimo komanda gale!
\address[B]{Department of Statistics\\
University of Toronto\\
100 St. George St., Rm. 6018\\
Toronto, Ontario\\
Canada M5S 3G3\\
\printead{u2}}
\end{aug}

% HISTORY:
\received{\smonth{1} \syear{2011}}
\revised{\smonth{4} \syear{2012}}

% ABSTRACT
%
\begin{abstract}
We analyze the convergence rate of a simplified version of a popular
Gibbs sampling method used for statistical discovery of gene
regulatory binding motifs in DNA sequences. This sampler satisfies a
very strong form of ergodicity (uniform). However, we show that, due
to multimodality of the posterior distribution, the rate of
convergence often decreases exponentially as a function of the length
of the DNA sequence. Specifically, we show that this occurs whenever
there is more than one true repeating pattern in the data. In practice
there are typically multiple such patterns in biological data, the
goal being to detect the most well-conserved and frequently-occurring
of these. Our findings match empirical results, in which the
motif-discovery Gibbs sampler has exhibited such poor convergence that
it is used only for finding modes of the posterior distribution
(candidate motifs) rather than for obtaining samples from that
distribution. Ours are some of the first meaningful bounds on the
convergence rate of a Markov chain method for sampling from a
multimodal posterior distribution, as a function of statistical
quantities like the number of observations.
\end{abstract}

% KEYWORDS
\begin{keyword}[class=AMS]
\kwd[Primary ]{62F15}
\kwd[; secondary ]{60J10}
\end{keyword}

\begin{keyword}
\kwd{Gibbs sampler}
\kwd{DNA}
\kwd{slow mixing}
\kwd{spectral gap}
\kwd{multimodal}
\end{keyword}

\end{frontmatter}

%s1 #&#
\section{Introduction}

\begin{figure}

\includegraphics{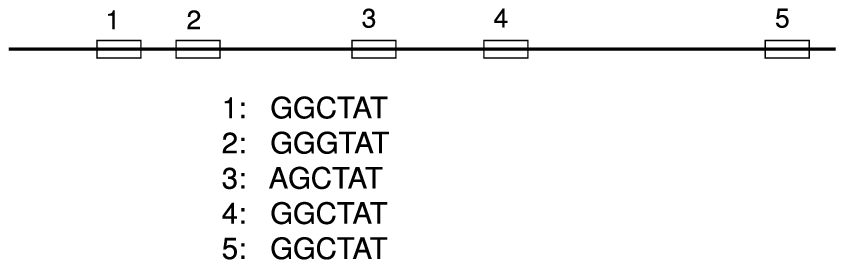}

\caption{Illustration of motif discovery: finding an unknown repeating
pattern in a long DNA sequence. The pattern can vary
slightly between instances.}\label{FigGoal}
\end{figure}

Gene regulatory binding motifs are short DNA sequences that control
gene expression. The identification of these regulatory motifs poses
several challenges: they are only 6--15 base pairs in length, and do
not contain clear start and stop codons; a regulatory motif is
indistinguishable from random sequences of the same length except that
it is a particular sequence that occurs more frequently than expected
under the background model. Discovery of previously undescribed
regulatory motifs in DNA sequences thus involves both finding such a
repeating pattern (``motif'') and determining where that pattern
occurs in the sequences [Kellis et al. (\citeyear{kellpattbirr04})];
this is illustrated in
Figure~\ref{FigGoal}.

%f1 #&#

One of the most effective methods for identifying new
regulatory motifs is based on a statistical model and associated Gibbs
sampling computational method [\citet{liuneuwlawr1995}]. This
approach has been popularized with the availability of software
programs for its use, such as BioProspector [\citet{liubrutliu01}] and
AlignAce [\citet{rothhughestechur98}].

Like most other methods for identifying regulatory motifs, the Gibbs
sampling method often yields different answers when starting from
different initial configurations. The method is applied by rerunning
the Gibbs sampler many times, using randomly generated initial
positions. The resulting candidate motifs are sorted according to
some goodness-of-fit measure, and then the highest-scoring motifs are
reported [\citet{lawraltsbogu93},
\citet{liubrutliu01},
\citet{jensliuzhou2004}].
This fact contrasts with the theoretical properties and traditional
use of a Gibbs sampler, namely to be simulated until it has some claim
of having converged to the posterior distribution, at which point the
answer should be the same regardless of initialization.

\begin{figure}

\includegraphics{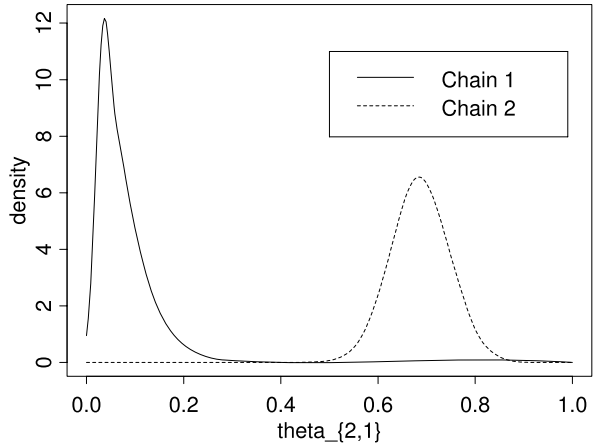}

\caption{The posterior density estimates of $\hat{\theta}_{2,1}(\bA
)$ from two different
Gibbs sampling chains, in the case of two true motifs.} \label{FigDensPlot}
\end{figure}

We address a particular model and Gibbs sampler that are
representative of this family of methods. We analyze the convergence
rate of a simplified version of the Gibbs sampler and show that, due
to multimodality of the posterior distribution, the convergence rate
typically decreases exponentially as a function of the DNA sequence
length (Theorem~\ref{ThmSlowMix2}). Specifically this occurs when
there is more than one true repeating pattern in the data, meaning
that the DNA is made up of short subsequences, each of which is either
equal to one of several motifs or is generated from the background
model. In practice there are typically multiple distinct repeating
patterns in biological data, corresponding to multiple gene regulatory
binding motifs or to repeating patterns that have other biological
significance, such as ``determinants of mRNA stability or even sites
for regulation by antisense transcripts''
[\citet{rothhughestechur98}]. The goal is to detect the most
frequently-occurring and well-conserved motif or motifs
[\citet{neuwliulawr95}]. So in practice we can expect the sampler
convergence rate to decay exponentially; this is equivalent to the run
time of the algorithm growing exponentially in the sequence length,
for a fixed accuracy. The multimodality of the posterior and
resulting poor convergence are illustrated in
Figure~\ref{FigDensPlot}, which shows posterior density estimates of
a particular function of the parameter vector, from two different
Gibbs sampling chains. Initialized with distinct parameter values,
the two chains have become trapped in different modes of the posterior
density and thus have not yet individually converged to the posterior
distribution.

%f2 #&#

The multimodality of the posterior distribution arises due to a
contradiction between the data, which typically have multiple true
repeating patterns, and the model assumption of a single such pattern.
Practitioners use the model not because it is believed to precisely
capture the true process that generated the data (which is extremely
complex) but because it captures several important features of that
process [\citet{neuwliulawr95},
\citet{rothhughestechur98}]. Our results
show that the presence of multiple motifs, even if some occur very
infrequently, causes slow convergence. Recognizing that there can be
multiple true motifs, a variant on the Gibbs sampler has been proposed
that allows for a fixed number of motifs greater than one
[\citet{neuwliulawr95}]. This approach is only likely to fix the
slow convergence if the number of motifs in the model is at least as
large as the number of true motifs in the data. This is only a
practical solution if the number of true motifs is small.

Our simplification of the model and associated Gibbs sampler assumes
that motifs can only end at locations in the sequence that are
divisible by the motif length, instead of at arbitrary locations
(Section~\ref{SecSimplification}). This is done to facilitate
analysis, by avoiding the ``phase shift'' issue that occurs in the
original sampler [\citet{lawraltsbogu93},
\citet{liu1994}]. Since
phase shift slows convergence of the chain, it is likely (but
unproven) that our results on slow convergence of the simplified chain
also hold for the original chain.

We also give evidence supporting the conjecture that the convergence
rate decreases polynomially if there is no more than one true (and
identifiable) motif in the data. We give empirical evidence, and
prove polynomial decay of the convergence rate for the case of
length-one motifs. In this case any true motifs are
nonidentifiable; see Theorem~\ref{ThmRapidMix}.

Ours are some of the few meaningful bounds on the convergence rate of
a Markov chain method used in Bayesian statistics, as a function of
statistical quantities such as the number of observations or number of
groups. Such results are particularly rare for multimodal posterior
densities. \citet{robesahu01} show that the convergence rate of a
Gibbs sampler for a unimodal posterior density in $\R^d$ approaches a
constant as the number of observations increases. \citet{bellcher09}
show that if the posterior density converges uniformly to a normal
density, then a Metropolis--Hastings chain restricted to a neighborhood
of the true parameter value has polynomially decaying convergence
rate. Jones and Hobert (\citeyear{Jones2001,jonehobe04}) and
other authors [e.g., Rosenthal (\citeyear{rose1995b,rose96})]
obtain bounds on the time to be within distance $\varepsilon> 0$ of
convergence for various hierarchical random effect models having
unimodal posterior densities, as a function of the initial values,
data and hyperparameters. \citet{mossvigo06} show that the
convergence rate of a Markov chain method used in Bayesian
phylogenetics can decrease exponentially in the number of samples in
the dataset. We also learned after completing this article that Dr.
Scott Schmidler at Duke University has independently obtained some
convergence results in the motif-discovery context (personal
communication).

Showing that a Markov chain method used in statistical practice is
``well-behaved'' usually consists of proving geometric ergodicity
[\citet{liuwongkong95},
\citet{jarnhans00},
\citet{fortmoulrobe03},
\citet{johnjone10}], meaning that the chain converges to the
posterior distribution at a geometric rate. The Gibbs sampler we
analyze satisfies the even stronger property of uniform
ergodicity; despite this, it is so poorly-behaved as to be unusable
for obtaining samples from the posterior distribution for long DNA sequences.

Characterizing the dependence of the convergence rate on statistical
quantities like the number of observations or the number of parameters
is critical in justifying the use of a Markov chain method. However,
there are several difficulties in doing so. First, the posterior
distribution of a statistical model has a much more complex form than
the stylized distributions for which Markov chain convergence rates
are typically obtained [Borgs et al. (\citeyear{borgchayfrie1999}),
Bhatnagar and Randall
(\citeyear{bhatrand2004}),
Woodard, Schmidler and Huber (\citeyear{woodschmhube20073})].
Second, the data, and thus the convergence rate of the Markov
chain, are stochastic and depend on the data-generating mechanism.

We address these challenges by utilizing Bayesian asymptotic theory,
which characterizes the behavior of the posterior distribution as the
number of observations grows. This is complicated by the fact that
Bayesian asymptotic theory is most well developed in the case of a
continuous parameter space, but the motif Gibbs sampler is defined on
a discrete parameter space. We handle this by applying the asymptotic
results on an alternative continuous parameterization of the motif
model and then mapping those results to the discrete parameterization.
Due to these technical challenges our main theorem requires
sufficiently long motifs, and is restricted to the case where each true
motif corresponds to a fixed sequence of nucleotides (disallowing the
small variations seen in Figure~\ref{FigGoal}). We give an additional argument
and simulation results suggesting that slow mixing holds even for very
short motifs, and when the true motifs are allowed to vary between
instances.

The motif discovery example provides insights into the dynamics of
standard Markov
chain methods applied to statistical models with highly multimodal
posterior distributions. Other examples that may have the same
exponential-time property include Markov chains for model search in
the context of regression with a large number of predictors
[\citet{lianwong2000},
\citet{hansdobrwest07}] and Markov chains for spatial mixture
models based on random fields
[\citet{gemagema1984},
\citet{greerich2002}]. Our example also provides a test case for the use
of more sophisticated Markov chain methods that are designed to handle
multimodality [\citet{delmdoucjasr06},
\citet{andrdouchole10}]: if a method can be shown to
sample from the posterior distribution of the motif-discovery model in
polynomial time, then it would be dramatically more efficient than the
Gibbs sampling approach.

Background on the Gibbs sampling method for motif discovery and on
Markov chain convergence rates is in Section~\ref{SecBackground}.
Our convergence results are in Section~\ref{SecMixRes}, and a
simulation study is given in Section~\ref{SecSimul}. The proof of
our main result is in Section~\ref{SecModeThm}, and we draw
conclusions in Section~\ref{SecConclude}.

%s2 #&#
\section{Background} \label{SecBackground}
%s2.1 #&#
\subsection{Statistical motif discovery}\label{SecMotifFinding}

The goal of motif discovery is to find short sub-sequences
of nucleotides (length 6--15 base pairs) that occur multiple times
(more often than could be explained under the background model) in one or
more long DNA sequences. Neither the nucleotide pattern nor
the sub-sequence locations are known. This goal is illustrated in
Figure~\ref{FigGoal}.

We address one of the two main variants of Gibbs sampler used in motif
discovery. The variant we analyze takes the number of motif
instances per sequence to be unknown, while the other variant fixes
the number of instances per sequence [\citet{jensliuzhou2004}]; the
two approaches are closely related and should have similar properties.
Programs such as BioProspector are based on the method we analyze, and
build in a number of additional features, such as a prior
distribution on the motif frequency and handling of gapped motifs;
however, by adding parameters and complexity to the model these
enhancements probably make the Gibbs sampler slower to converge,
and so are unlikely to affect our slow-mixing results.

We focus further on the case of a single DNA sequence (having an unknown
number of motif instances). The case of multiple sequences can be
addressed by concatenating to obtain a single sequence.

The motif instances are not necessarily identical. Taking the length
$w$ of the motif to be known, one can describe the nucleotide pattern
by a position-specific frequency matrix, which contains the
probability of occurrence of each nucleotide at each position in the
motif. Call this matrix $\btheta_{1:w} = (\btheta_1, \ldots, \btheta_w)$,
where $\btheta_k$ is the unknown probability vector for the $k$th
position. Let the nucleotides be labeled $1,\ldots,M$, so that
$\btheta_k$ has length $M$; for DNA data $M=4$. For each instance of
the motif, the nucleotide in position $k$ is assumed to be drawn
independently from a discrete distribution with parameters
$\btheta_k$. The positions in the full sequence that are not part of
a motif instance are assumed to have nucleotide drawn independently
from a discrete distribution with unknown probability vector
$\btheta_0$.

%We are considering the case of multiple DNA sequences where the motif
%occurs an unknown number of times in each sequence. This can be
%simplified to the case of a single sequence by concatenating the
%separate sequences. This ignores edge effects but will not affect the
%dynamics of the Markov chain we will define, when the sequences are
%very long.

Let $\bS= (S_1, \ldots, S_{L}) \in\{1,\ldots,M\}^L$ be the observed
sequence, having\break length~$L$. In the original model of, for example,
\citet{liuneuwlawr1995}, a motif is allowed to start at any index
$i\in\{1, \ldots, L-w+1\}$, but we will analyze a simplified version
that only allows a motif to start at indices $wi-w+1$ for $i \in\{1,
\ldots, L/w\}$ where $L$ is divisible by $w$. This choice is explained
in Section~\ref{SecSimplification}. Let $A_i \in\{0,1\}$ be the (unknown)
indicator of whether a motif begins at index $wi-w+1$, for $i\in\{1,
\ldots, L/w\}$, and define $\bA= (A_1, \ldots, A_{L/w})$. Let
$\bN(\bA^{(k)})$ be the length-$M$ vector of counts of the occurrence
of each nucleotide at position $k \in\{1,\ldots,w\}$ of all motif
instances, conditional on $\bA$. Similarly, $\bN(\bA^c)$ is defined to
be the length-$M$ vector of counts for each nucleotide in the
background locations, that is, the locations that are not part of any
motif instance,
%
%e2.1 #&#
\begin{eqnarray}
\label{EqnCountVecs}\qquad N \bigl(\bA^{(k)} \bigr)_m &
\triangleq& \sum_{i=1}^{L/w}
\indic_{\{A_i=1,
S_{wi-w+k}=m\}},\qquad k \in \{1,\ldots,w\}, m\in\{1,\ldots, M\},
\nonumber\\
N \bigl(\bA^c \bigr)_m &\triangleq& N(
\bS)_m - \sum_{k=1}^w N \bigl(
\bA^{(k)} \bigr)_m,
\\
N(\bS)_m &\triangleq& \sum_{i=1}^L
\indic_{\{S_{i}=m\}}.
\nonumber
\end{eqnarray}

For any two equal-length vectors $\bbeta= (\beta_1, \ldots, \beta_K)$ and $\mathbf{n}
= (n_1, \ldots, n_K)$ define the notation
%
%e2.2 #&#
\begin{equation}
\label{EqnNotation} \bbeta^{\mathbf{n}} = \prod_{k=1}^K
\beta_k^{n_k},\qquad \Gamma(\mathbf{n}) = \prod
_{k=1}^K \Gamma(n_k),\qquad |\mathbf{n}| =
\sum_{k=1}^K n_k,
\end{equation}
where $\Gamma$ is the gamma function. Using this notation, the
likelihood conditional on $\bA$ can be written as
%
%e2.3 #&#
\begin{eqnarray}\label{EqnFullLike}
\pi(\bS| \bA, \btheta_{0:w}) &= &\prod_{i=1}^{L/w}
\Biggl[ \prod_{k=1}^w
\theta_{k,S_{wi-w+k}} \Biggr]^{\indic_{\{A_i=1\}}} \Biggl[ \prod
_{k=1}^w \theta_{0,S_{wi-w+k}} \Biggr]^{\indic_{\{A_i=0\}}}
\nonumber
\\[-8pt]
\\[-8pt]
\nonumber
&= &\btheta_0^{\bN(\bA^c)} \times \prod
_{k=1}^w \btheta_k^{\bN(\bA^{(k)})},
\nonumber
\end{eqnarray}
where $\btheta_{0:w} = (\btheta_0, \btheta_1, \ldots, \btheta_w)$
denotes all model parameters. We will use $\pi$ to indicate the
likelihood, prior or the full, marginal or conditional posterior
distributions as distinguished by its arguments.

The prior distributions for the unknown quantities are $\pi(\btheta_k)
= \operatorname{Dirichlet}(\bbeta_k)$ for $k \in\{0,\ldots,w\}$ and $\pi(A_i
= 1) = p_0$ independently. Here $p_0 \in(0,1)$ is a known constant
and $\bbeta_0, \ldots, \bbeta_w$ are fixed length-$M$ vectors with
$\beta_{k,m} > 0$. The corresponding posterior distribution is
[\citet{jensliuzhou2004}]
%
%e2.4 #&#
\begin{eqnarray}\label{EqnFullPost}
&&\pi(\bA, \btheta_{0:w} | \bS) \propto\pi(\bA, \btheta_{0:w} ,
\bS)
\nonumber
\\[-8pt]
\\[-8pt]
\nonumber
&&\qquad= p_0^{|\bA|} (1-p_0)^{L/w-|\bA|} \times
\btheta_0^{\bN(\bA^c)+\bbeta_0-1} \times\prod_{k=1}^w
\btheta_k^{\bN(\bA^{(k)}) + \bbeta_k - 1}.
\nonumber
\end{eqnarray}

\citet{liuneuwlawr1995} integrate out the parameters
$\btheta_{0:w}$ in the above formula to yield a posterior
distribution on $\bA$. Using the notation from (\ref{EqnNotation}),
%
%e2.5 #&#
\begin{eqnarray}\label{EqnPostDistn}
\qquad &&\pi(\bA| \bS)
\nonumber
\\[-8pt]
\\[-8pt]
\nonumber
&&\qquad\propto p_0^{|\bA|} (1-p_0)^{L/w-|\bA|}
\frac{\Gamma( \bN(\bA^c) +
\bbeta_0)}{\Gamma(|\bN(\bA^c)|+|\bbeta_0|)} \prod_{k=1}^w
\frac{\Gamma(\bN(\bA^{(k)}) +
\bbeta_k)}{\Gamma(|\bN(\bA^{(k)})| + |\bbeta_k|)}.
\nonumber
\end{eqnarray}
\citet{liu1994} gives theoretical results supporting faster
convergence of a Gibbs sampler for the reduced posterior
$\pi(\bA|\bS)$ relative to a Gibbs sampler for
$\pi(\bA,\btheta_{0:w}|\bS)$. So \citet{liuneuwlawr1995}
propose to use a Gibbs sampler to draw from $\pi(\bA|\bS)$, having
state vector $\bA\in\cX\triangleq\{0,1\}^{L/w}$. This sampler
iteratively updates each $A_i$ according to its conditional posterior
distribution; details are given in Section~\ref{SecGibbs}.

Although this Gibbs sampler has both systematic-scan and random-scan
versions, we expect that the mixing properties (defined in
Section~\ref{SecMixingTime}) of the two versions are identical. For
this reason we focus attention on the random-scan Gibbs sampler, which
is easier to analyze. We also make the transition matrix of the
Markov chain nonnegative definite by including a holding probability
of $1/2$ at every state; this is a common technique when analyzing the
mixing properties of Markov chains
[\citet{madrzhen2003},
\citet{woodschmhube09a}]. It only increases the
mixing time (Section~\ref{SecMixingTime}) by a factor of two, so it does
not affect results on the order of the run time as a function of $L$.
Let $\bA_{[-i]}$ indicate the vector $\bA$ excluding the $i$th
element, and $\pi(A_i | \bA_{[-i]}, \bS) \propto\pi(\bA| \bS)$
indicate the conditional posterior distribution of $A_i$ given
$\bA_{[-i]}$. With these definitions we can write the transition
matrix $T$ of the Gibbs sampler as follows for any $\bA, \bA' \in
\cX$:
\begin{eqnarray}
\label{EqnTransMat} T \bigl(\bA, \bA' \bigr) &\triangleq&
\frac{1}{L/w} \sum_{i=1}^{L/w}
\indic_{\{\bA'_{[-i]} =
\bA_{[-i]}\}}
\nonumber\\
&&\hspace*{41pt}{} \times\frac{1}{2} \bigl[\indic_{\{A_i' = A_i\}} + \pi(A_i = 0 |
\bA_{[-i]}, \bS) \indic_{\{A_i'=0\}}\\
&&\hspace*{102pt}{} + \pi(A_i = 1 |
\bA_{[-i]}, \bS) \indic_{\{A_i'=1\}} \bigr].\nonumber
\end{eqnarray}
Expressions for $\pi(A_i | \bA_{[-i]}, \bS)$ are given in
Section~\ref{SecGibbs}.

We also need an expression for the likelihood marginalized over
$\bA$. Using the vector notation $\bS_{n:m} = (S_n,\ldots,S_m)$ for
$n\leq m$,
%
%e2.6 #&#
%e2.7 #&#
\begin{eqnarray}\quad
\pi(\bS| \btheta_{0:w}) %&=
% \prod_{i=1}^{L/w} \sum_{A_i =0}^1 \left[ p_0 \prod_{k=1}^w
% \theta_{k,S_{wi-w+k}}\right]^{\indic_{\{A_i=1\}}} \left[ (1-p_0)
&=& \prod_{i=1}^{L/w} f(
\bS_{wi-w+1:wi} | \btheta_{0:w}) \qquad\mbox {where}\label{EqnMargLike}
\\
f(\mathbf{s}| \btheta_{0:w}) &\triangleq& p_0 \prod
_{k=1}^w \theta_{k,s_k} +
(1-p_0) \prod_{k=1}^w
\theta_{0,s_k},\qquad \mathbf{s}\in\{1,\ldots, M\}^w. \label{EqnFSDef}
\end{eqnarray}

Under model (\ref{EqnMargLike}) and (\ref{EqnFSDef}) each subsequence
$\bS_{wi-w+1:wi}$ is either generated from the motif $\btheta_{1:w}$ with
probability $p_0$, or generated from the background $\btheta_0$ with
probability $(1-p_0)$. So we have i.i.d. observations
$\bS_{wi-w+1:wi}$ for $i \in\{1,\ldots,L/w\}$, which will allow us to
use Bayesian asymptotic theory for i.i.d. parametric models.

%s2.2 #&#
\subsection{Reason for the simplification}\label{SecSimplification}
As stated in Section~\ref{SecMotifFinding}, while the original model
allows a motif to start at any index $i \in\{1, \ldots, L-w+1\}$, we
analyze a simplification of the model that assumes that motifs can
only start at indices $wi-w+1$ for $i\in\{1, \ldots, L/w\}$. This
simplification is done to facilitate analysis; however, we believe
that our results are likely to hold for the original model as well.

First, our rapid mixing result (Theorem~\ref{ThmRapidMix})
immediately holds for the original model and associated Gibbs sampler.
This is because it is for the case of $w=1$, where the original and
simplified models are identical.

Second, the proof of our slow mixing result
(Theorem~\ref{ThmSlowMix2}) can be extended to the case where the
model allows motifs to start at indices that are fixed distance $ \geq
w$ apart. However, Theorem~\ref{ThmSlowMix2} does not easily
extend to the case where motifs can start in locations that are less
than $w$ distance apart (including the original model). This is due
to the following ``phase shift'' issue, which
complicates analysis. For illustration consider the case where $M=4$
(there are four possible nucleotides) and $w=5$ (motifs are five
nucleotides long), and the true motif is (deterministically) the
sequence $(1,4,2,2,3)$. Phase shift means that it is possible to
estimate that a motif begins or ends in the middle of one of the
$(1,4,2,2,3)$ subsequences that exist in the data. For example, if
the DNA sequence $\bS$ satisfies $\bS_{22:26} = (1,4,2,2,3)$,
corresponding to a true motif beginning at position $22$, then the
original model also allows for the possibility that $A_{23} = 1$,
meaning that a motif could be estimated to instead start at position
$23$ with the sequence $4,2,2,3$.

While phase shift complicates analysis of the original Gibbs sampler,
it should also make the original Gibbs sampler converge more slowly
than the simplified Gibbs sampler. The effect of phase shift on
convergence of the original Gibbs sampler is that it can become trapped in
a local mode of the posterior distribution that corresponds to a
shifted version of the true motif. This effect is described in
\citet{lawraltsbogu93} and \citet{liu1994}. To
illustrate, take the above example where the true motif is
$(1,4,2,2,3)$. There is a local mode of the
posterior distribution for which the inferred motif starts with the
sequence $4,2,2,3$, another for which the inferred motif ends with
the sequence $1,4,2,2$, and so on. This posterior multimodality slows
convergence of the
original Gibbs sampler. It also suggests that our slow mixing result
(Theorem~\ref{ThmSlowMix2}) for the simplified Gibbs sampler
holds for the original sampler.

Analysis of the original sampler should be possible
using the same general approach taken here, but a number of the
technical details would need to change. We leave this to future work.

%s2.3 #&#
\subsection{Markov chain convergence rates}\label{SecMixingTime}

Consider a Markov chain with transition matrix $T$ and stationary
distribution $\pi$ on a discrete state space $\cX$. For $x \in\cX$
and $D \subset\cX$, let $T(x, D) = \sum_{y \in D} T(x,y)$. If the
chain is initialized at $x \in\cX$, then the total variation distance
to stationarity after $n$ iterations is
\[
\bigl\|T^n(x,\cdot) - \pi(\cdot)\bigr\|_{\mathrm{TV}} \triangleq
\max_{D \subset\cX}\bigl | T^n(x,D)-\pi(D) \bigr|.
\]
The mixing time of the chain is the number of iterations required to
be within distance $\varepsilon\in(0,1)$ of stationarity,
\begin{eqnarray*}
\tau_\varepsilon\triangleq\max_{x \in\cX} \min \bigl\{n\dvtx
\bigl\|T^{m}(x,\cdot ) - \pi(\cdot)\bigr\|_{\mathrm{TV}} \leq\varepsilon\mbox{
for all } m\geq n \bigr\};
\end{eqnarray*}
cf. Sinclair (\citeyear{sinc1992}).\vadjust{\goodbreak}

Consider $T$ irreducible, aperiodic,
reversible and nonnegative definite, which holds for the random-scan Gibbs
sampler in Section~\ref{SecMotifFinding}. Then $\tau_\varepsilon$ is
finite and
closely related to the spectral gap $\Gap(T) \triangleq1 - \lambda_2,$
where $\lambda_2 \in[0,1)$ is the second-largest eigenvalue of $T$.
Since the state space $\cX$ is finite, $\Gap(T) > 0$ and the chain is
called \textit{uniformly ergodic} [\citet{roberose2004}]. The quantities
$\tau_\varepsilon$ and $\Gap(T)$ are related via [\citet{sinc1992}]
%
%e2.8 #&#
\begin{eqnarray}\label{EqnMixTimeGap}
\tau_\varepsilon&\leq&\Gap(T)^{-1} \Bigl( -\ln \Bigl[
\min_{x\in\cX} \pi(x) \Bigr] - \ln\varepsilon \Bigr),
\nonumber
\\[-9pt]
\\[-9pt]
\nonumber
\tau_\varepsilon&\geq&\frac{1}{2} \bigl(1-\Gap(T) \bigr)
\Gap(T)^{-1} \bigl( - \ln (2 \varepsilon) \bigr).
\end{eqnarray}

The efficiency of the Markov chain can be measured by how quickly
$\tau_\varepsilon$ increases as a function of the problem difficulty, for
instance the dimension of the parameter space. In our case we are
interested in the dependence of $\tau_\varepsilon$ on the length~$L$ of
the DNA sequence, since in practice one analyzes very long sequences.
We would certainly hope that $\tau_\varepsilon$ grows at most
polynomially in $L$ for any fixed $\varepsilon$; this property is called
\textit{rapid mixing}. \textit{Slow mixing} means that $\tau_\varepsilon$
increases exponentially for some $\varepsilon$. By
(\ref{EqnMixTimeGap}) rapid mixing is equivalent to $\Gap(T)$
decreasing at most polynomially toward zero, and slow mixing is
equivalent to $\Gap(T)$ decreasing exponentially toward zero, if
$-\log[\min_{x\in\cX} \pi(x)]$ increases polynomially in~$L$. The
latter property holds for the random-scan Gibbs sampler in
Section~\ref{SecMotifFinding}. The rapid/slow mixing distinction is a
measure of the computational tractability of an algorithm; polynomial
factors are expected to be eventually dominated by increases in
computing power due to Moore's Law, while exponential factors cause a
persistent computational problem.\vspace*{-2pt}

%s3 #&#
\section{Convergence results}\label{SecMixRes}

We consider the mixing time (equivalently, spectral gap) of the Gibbs
sampler when the data are drawn from a generalization of the model
given in
Section~\ref{SecMotifFinding} that allows multiple true motifs. First
we give negative results for the case of multiple true motifs, and then we
give a positive result for a case with no true motifs.\vspace*{-2pt}

%s3.1 #&#
\subsection{Slow mixing for multiple true motifs} \label{SecMixRes1}

In this section, we show that if the data actually contain multiple
true motifs, then the Gibbs sampler is slowly mixing: that $\Gap(T)$
is $\mathcal{O}(\alpha^{-L})$ for $\alpha> 1$ where $-\log[\min_{x
\in\cX} \pi(x)]$ is $\mathcal{O}(L^q)$ for some $q > 0$. To make
this statement precise in the presence of random data, we need to make
assumptions about the model by which the data are generated. Our
convergence results are obtained using Bayesian asymptotics based on
this generative model; for other connections between Markov chain
convergence and Bayesian asymptotics, see \citet{kama11} and
\citet{bellcher09}.

For a concrete example to keep in mind, consider the case where $M=4$
(there are four possible nucleotides)\vadjust{\goodbreak} and $w=5$ (motifs are five
nucleotides long). Then let the DNA sequence $\bS$ be generated as
the concatenation of many length-five subsequences, each of which is
either (motif one) equal to $(1,4,2,2,3)$ with probability 0.005, or
(motif two) equal to $(4,2,4,1,3)$ with probability 0.001, or
generated as i.i.d. noise where each nucleotide has equal
probability. Theorem~\ref{ThmSlowMix2} below says that the Gibbs
sampler $T$ is slowly mixing for data generated in this way.

When analyzing the Gibbs sampler we do \textit{not} assume that the data
$\bS$ are generated according to the inference model
(\ref{EqnMargLike}) and (\ref{EqnFSDef}). Our most general result only
assumes that the subsequences $\bS_{wi-w+1:wi}$ are i.i.d.
%
%as3.1 #&#
\begin{assumption}\label{AssuGener}
The data subsequences $\bS_{wi-w+1:wi}$ indexed by $i \in
\{1,\ldots, L/w\}$ are independent and identically distributed
according to some probability mass function
$g(\mathbf{s})>0\dvtx \mathbf{s}\in\{1,\ldots,M\}^w$, that is,
\[
\bS\sim\prod_{i=1}^{L/w} g(
\bS_{wi-w+1:wi}).
\]
\end{assumption}
Under Assumption~\ref{AssuGener}, we give a simple sufficient
condition for slow mixing that relates the generative model $g(\mathbf{s})$
to the inference model $f(\mathbf{s}| \btheta_{0:w})$ via the
quantity $E
\log f(\mathbf{s}|\btheta_{0:w}) = \sum_{\mathbf{s}} g(\mathbf{s})
\log f(\mathbf{s}|
\btheta_{0:w})$. Since $\btheta_k$ for $k \in\{0,\ldots,w\}$ is
defined on the simplex $\Psi\triangleq\{\btheta_k\dvtx \sum_{m=1}^M
\theta_{k,m}=1, \theta_{k,m} \geq0\}$, the quantity $E \log
f(\mathbf{s}
| \btheta_{0:w}) \in[-\infty,0)$ is a function of $\btheta_{0:w}
\in
\Psi^{w+1}$. It is continuous in $\btheta_{0:w}$, because it is a
linear combination of a finite number of the continuous functions
$\log f(\mathbf{s}|\btheta_{0:w})$. We call $\eta(\btheta_{0:w})
\triangleq
E \log f(\mathbf{s}| \btheta_{0:w})$ \textit{multimodal} if there exist
$\btheta^{(1)}_{0:w},\btheta^{(2)}_{0:w} \in\Psi^{w+1}$ and bounded
sets $F_1
\ni\btheta^{(1)}_{0:w}$ and $F_2 \ni\btheta^{(2)}_{0:w}$ such that
%
%e3.1 #&#
\begin{equation}
F_1 \cap F_2 = \varnothing \quad\mbox{and}\quad
\label{EqnMultimodeDef} \sup_{\btheta_{0:w} \in\partial F_j} \eta(\btheta_{0:w}) < \eta
\bigl(\btheta_{0:w}^{(j)} \bigr), \qquad j \in\{1,2\},
\end{equation}
where $\partial F_j \triangleq\operatorname{cl}(F_j) \cap\operatorname{cl}(F_j^c)$
is the boundary of $F_j$. Equation (\ref{EqnMultimodeDef}) implies
that $\btheta^{(j)}_{0:w}$ is in the interior of $F_j$. For a
continuous function on a closed, connected subset of $\R^d$ (like
$\Psi^{w+1}$) this definition of multimodality is weaker than the
existence of multiple strict local maxima and stronger than the
existence of multiple local maxima. We call a function $h$ of
$\btheta_{0:w}$ a \textit{multiminimum function} if $-h(\btheta_{0:w})$
is multimodal.
%
%th3.1 #&#
\begin{theorem}\label{ThmSlowMix}
Under Assumption~\ref{AssuGener}, if the function $E \log f(\mathbf{s}|
\btheta_{0:w})$ of $\btheta_{0:w}$ is multimodal then the spectral gap
of the Gibbs sampler $T$ decreases exponentially in $L$, almost
surely.
\end{theorem}
Theorem~\ref{ThmSlowMix} will be proven in Section~\ref{SecModeThm}.
It uses asymptotic results on the behavior of the posterior when the
data are not distributed according to the inference model
[\citet{berk65}]. We will see that when $E \log f(\mathbf{s}|
\btheta_{0:w})$ is multimodal\vadjust{\goodbreak} the posterior distribution is also
multimodal for large $L$, and that the heights of the modes relative
to the heights of the valleys in between grow exponentially in $L$,
causing the slow mixing. This is due to the fact that, using
(\ref{EqnMargLike}), the log-likelihood is
\[
\log\pi(\bS| \btheta_{0:w}) = \sum_{i=1}^{L/w}
\log f(\bS_{wi-w+1:wi} | \btheta_{0:w}),
\]
which satisfies the following for any $\btheta_{0:w}$:
\[
\frac{1}{L/w} \sum_{i=1}^{L/w} \log f(
\bS_{wi-w+1:wi}| \btheta_{0:w}) \stackrel{L \rightarrow\infty} {
\longrightarrow} E \log f(\mathbf{s}| \btheta_{0:w}) \qquad\mbox{a.s.}
\]
by the strong law of large numbers. This effect leads to the likelihood
function being multimodal for large $L$. Statistically,
these correspond to multiple values of
$\btheta_{0:w}$ that explain the data well.

Another way of stating Theorem~\ref{ThmSlowMix} is via the
Kullback--Leibler divergence between $f(\mathbf{s}| \btheta_{0:w})$ and
$g(\mathbf{s})$. The divergence measures the degree of difference between
$f(\mathbf{s}| \btheta_{0:w})$ and $g(\mathbf{s})$ and is defined as
%
%e3.2 #&#
\begin{eqnarray}\label{EqnKLdiverge}
&&\sum_{\mathbf{s}} g(\mathbf{s}) \log\frac{g(\mathbf
{s})}{f(\mathbf{s}| \btheta_{0:w})}
\nonumber
\\[-8pt]
\\[-8pt]
\nonumber
&&\qquad= \sum_\mathbf{s}g(\mathbf{s}) \log g(\mathbf{s}) -
\sum_{\mathbf
{s}} g(\mathbf{s}) \log f(\mathbf{s}|
\btheta_{0:w}).
\end{eqnarray}
Since $\sum_\mathbf{s}g(\mathbf{s}) \log g(\mathbf{s})$ does not
depend on
$\btheta_{0:w}$, the divergence is a multiminimum function iff $E \log
f(\mathbf{s}| \btheta_{0:w})$ is multimodal.
%
%co3.1 #&#
\begin{cor} \label{ThmKL}
Under Assumption~\ref{AssuGener}, if the Kullback--Leibler
divergence~(\ref{EqnKLdiverge}) is a multiminimum function of
$\btheta_{0:w}$ then the spectral gap of the Gibbs sampler $T$
decreases exponentially in $L$, almost surely.
\end{cor}

Next we show that multimodality of $E \log f(\mathbf{s}| \btheta_{0:w})$
occurs when the generative model $g(\mathbf{s})$ includes $J>1$ true motifs,
described by position-specific frequency matrices
$\btheta_{1:w}^{j*}$ for $j \in\{1,\ldots,J\}$.
Assumption~\ref{AssuMultiMotifs} says that $g(\mathbf{s})$ is
obtained by
extending the inference model (\ref{EqnFSDef}) to the case of $J>1$
motifs.

%as3.2 #&#
\begin{assumption} \label{AssuMultiMotifs}
The p.m.f. $g(\mathbf{s})$ from Assumption~\ref{AssuGener} satisfies
%
%e3.3 #&#
\begin{equation}
\label{EqnGSDef} g_{\btheta^*}(\mathbf{s}) = \sum
_{j=1}^J p_j \prod
_{k=1}^w \theta_{k,s_k}^{j*} +
\Biggl(1-\sum_{j=1}^J p_j
\Biggr) \prod_{k=1}^w \theta_{0,s_k}^*,
\end{equation}
where $J>1$ and:

\begin{longlist}[(1)]

\item[(1)] $p_j
>0$ for $j\in\{1,\ldots,J\}$ are the motif frequencies, where $\sum_{j=1}^J p_j < 1$;

\item[(2)] $\btheta_0^*$ is a background probability vector with
$\sum_{m=1}^M
\theta_{0,m}^* = 1$ and \mbox{$\theta_{0,m}^* > 0$};

\item[(3)] $\btheta_{1:w}^{j*}$ for $j \in\{1,\ldots,J\}$ are
position-specific frequency matrices.\vadjust{\goodbreak}
\end{longlist}

\end{assumption}

Due to the complex form of $E \log f(\mathbf{s}| \btheta_{0:w}) =
\sum_{\mathbf{s}}
g_{\btheta^*}(\mathbf{s}) \log f(\mathbf{s}| \btheta_{0:w})$ under
Assumptions~\ref{AssuGener} and~\ref{AssuMultiMotifs} it is difficult to
characterize the number of modes without making any additional
assumption. We will restrict our analysis to the case where for each $j
\in
\{1,\ldots,J\}$ and $k \in\{1,\ldots,w\}$ there is some $m \in
\{1,\ldots,M\}$ with $\theta_{k,m}^{j*} = 1$. This means that each
true motif is a fixed length-$w$ sequence of nucleotides, for example,
where $w=5$ and $M=4$ and the first and second motifs correspond to
the deterministic sequences $(1,4,2,2,3)$ and $(4,2,4,1,3)$,
respectively. The case without this restriction is discussed below.\vspace*{-2pt}

%as3.3 #&#
\begin{assumption}\label{AssuDeter}
For each $j \in\{1,\ldots,J\}$ and each $k \in
\{1,\ldots,w\}$, there is some $t^j_k \in\{1,\ldots,M\}$ for which
$\theta_{k,t^j_k}^{j*} = 1$. Also,
%
%e3.4 #&#
\begin{equation}
\label{EqnADef} a \triangleq\min_{j \neq j'} \liminf_{w\rightarrow\infty}
\frac{1}{w}\sum_{k=1}^w
\indic_{ \{t_k^j \neq t_k^{j'} \}} > 0.\vspace*{-2pt}
\end{equation}
\end{assumption}
Assumption~\ref{AssuDeter} says that the motifs are
deterministic in the above sense, and that for any two
motifs $j \neq j'$ the proportion of differences between
the motif sequences does not decay to zero as $w$ grows. This ensures
that the
motifs are different enough from one another for large $w$ to
%Define the quantity
%%
% j'} \sum_{k=1}^w \indic_{\{\btheta^{j*}_{k} \neq\btheta^{j'*}_k\}}
%%
%which is the smallest number of differences between any pair of true
%motifs.
%
%For our example with $J=2$ true motifs $(1,4,2,2,3)$ and
%$(4,2,4,1,3)$, we have $\mbox{diff}(\btheta^*) = 4$. We will give our
%slow mixing results for $\mbox{diff}(\btheta^*)$ large enough,
%ensuring that the motifs are different enough from one another to
cause a mixing problem in the Markov chain. With these assumptions,
we have multimodality of $E \log f(\mathbf{s}| \btheta_{0:w})$ for large
enough $w$.\vspace*{-2pt}
%
%le3.1 #&#
\begin{lemma}\label{ThmLocalMaxMotifs}
Under Assumptions~\ref{AssuGener}--\ref{AssuDeter} there exists
$w^*<\infty$ that depends on $p_0$, $\btheta_0^*$, $J$,
$\{p_j\}_{j=1}^J$, and $a$ [as in equation (\ref{EqnADef})] such that
the following holds. If $w \geq
w^*$, then $E \log f(\mathbf{s}| \btheta_{0:w})$ is multimodal with
at least
$J>1$ local maxima.\vspace*{-2pt}
\end{lemma}

Lemma~\ref{ThmLocalMaxMotifs} is proven in the supplementary material [\citet{woodrose12}]. Combining
Theorem~\ref{ThmSlowMix} and Lemma~\ref{ThmLocalMaxMotifs}
immediately yields our main result: slow mixing for the case of
multiple true motifs.\vspace*{-2pt}
%
%th3.2 #&#
\begin{theorem} \label{ThmSlowMix2}
Under Assumptions~\ref{AssuGener}--\ref{AssuDeter}, there
exists $w^* < \infty$ such that whenever $w \geq
w^*$ the spectral gap of the Gibbs sampler $T$ decreases exponentially
in $L$, almost surely.\vspace*{-2pt}
\end{theorem}
While Theorem~\ref{ThmSlowMix2} is stated for
$w$ large enough and assumes deterministic true motifs, the simulation results
in Section~\ref{SecSimul} suggest that slow mixing occurs even for
nondeterministic true motifs and for $w$ as low as six.

Theorem~\ref{ThmSlowMix2} says that the presence of multiple motifs
in the generative model contradicts the inference model assumption of
a single motif, causing slow mixing. In realistic biological
situations there are frequently multiple motifs, corresponding to
multiple gene regulatory binding motifs or to repeating patterns that
have other biological significance
[\citet{neuwliulawr95},
\citet{rothhughestechur98}].
Theorem~\ref{ThmSlowMix2} says that these\vadjust{\goodbreak} patterns do not have to
occur often in order to cause slow mixing (that slow mixing occurs
even when some of the $p_j$ are very small). So when $L$ is large the
Gibbs sampler should be used only as a tool for generating candidate
motifs, and the results cannot be interpreted as obtaining samples
from the posterior distribution, or used for Monte Carlo estimation.

Theorem~\ref{ThmSlowMix2} assumes that each true motif is a
deterministic sequence; now consider the case of variable true motifs,
that is, where Assumptions~\ref{AssuGener} and~\ref{AssuMultiMotifs} hold
but not Assumption~\ref{AssuDeter}. We give an informal argument
suggesting that the function $E \log f(\mathbf{s}| \btheta_{0:w})$ is still
multimodal.

Consider the case where $\sum_{j=1}^J
p_j = p_0$. Then, using (\ref{EqnFSDef}) and (\ref{EqnGSDef}),
\begin{eqnarray*}
g_{\btheta^*}(\mathbf{s}) &=& \sum_{j=1}^J
\frac{p_{j}}{p_0} \Biggl[ p_0 \prod_{k=1}^w
\theta_{k,s_k}^{j*} + (1-p_{0}) \prod
_{k=1}^w \theta_{0,s_k}^* \Biggr]
\\[-2pt]
&= &\sum_{j=1}^J \frac{p_j}{p_0} f
\bigl( \mathbf{s}| \bigl(\btheta_0^*, \btheta_{1:w}^{j*}
\bigr) \bigr).
\end{eqnarray*}
So $E \log f(\mathbf{s}| \btheta_{0:w})$ can be
written as
%
%e3.5 #&#
\begin{equation}\quad
\sum_\mathbf{s}g_{\btheta^*}(\mathbf{s}) \log f(
\mathbf{s}| \btheta_{0:w}) = \sum_{j=1}^J
\frac{p_{j}}{p_0} \sum_\mathbf{s}f \bigl(\mathbf{s}|
\bigl(\btheta_0^*, \btheta_{1:w}^{j*} \bigr)
\bigr) \log f(\mathbf{s}| \btheta_{0:w}) \label{EqnKLmixture}.
\end{equation}
By standard information-theoretic results [\citet{kull59},
\citet{berk65}], for each~$j$
the function $\sum_\mathbf{s}f(\mathbf{s}| (\btheta_0^*, \btheta_{1:w}^{j*}) ) \log f(\mathbf{s}|
\btheta_{0:w})$ has a unique global maximum at $\btheta_{0:w} =
(\btheta_0^*, \btheta_{1:w}^{j*})$.
Since (\ref{EqnKLmixture}) is the weighted sum of continuous
functions that have global maxima occurring at distinct locations
$(\btheta_0^*, \btheta_{1:w}^{j*})$, it seems likely that
(\ref{EqnKLmixture}) is multimodal when $J > 1$.\vspace*{-2pt}

%s3.2 #&#
\subsection{\texorpdfstring{Rapid mixing for $\leq$1 true motif}{Rapid mixing for <=1 true motif}}

Simulations suggest (Section~\ref{SecSimul}) that when there is no
more than one true motif the Gibbs sampler is rapidly mixing, that is,
$\Gap(T)^{-1} = \mathcal{O}(L^q)$ for some $q > 0$. We have one
theoretical result in this direction, showing rapid mixing for the
case $w=1$. In this case any true motif is indistinguishable from the
background signal, so there are effectively zero true motifs.\vspace*{-2pt}

%th3.3 #&#
\begin{theorem}\label{ThmRapidMix}
If $w=1$, then the Gibbs sampler $T$ has spectral gap that decreases
polynomially in $L$, uniformly over $\bS\in
\{1,\ldots,M\}^L$. Specifically, for
$M=2$,
\[
\sup_{\bS\in\{1,\ldots,M\}^L}\Gap(T)^{-1} = \mathcal{O} \bigl(L^{14}
\bigr)
\]
and for fixed $M>2$ the same result holds for a larger-degree polynomial.\vspace*{-2pt}
\end{theorem}
Theorem~\ref{ThmRapidMix} (proven in the supplementary material [\citet{woodrose12}]) shows
rapid mixing in the\vadjust{\goodbreak} worst
case over possible datasets $\bS$. Contrast with
Theorems~\ref{ThmSlowMix} and~\ref{ThmSlowMix2}, which show slow
mixing almost surely under a particular generative model $g(\mathbf{s})$.

It is likely that the spectral gap bound given in
Theorem~\ref{ThmRapidMix} is very loose as a function of $L$, since
the tools that we use to obtain it (Theorem~\ref{ThmPathBound} in
particular) can be imprecise. However, obtaining a tighter bound
would require substantially longer arguments, so we leave this to
future work. Additionally, one could assume a particular distribution for
$\bS$ and use an average-case analysis to obtain a tighter bound,
but also one that would have narrower interpretation.\vspace*{-2pt}

%s4 #&#
\section{Simulation study}\label{SecSimul}

We simulate data with either $J=1$ or $J=2$ true motifs, and measure
the convergence of the Gibbs sampler. The data are simulated as
follows; to emulate DNA data we take $M=4$. The true
position-specific frequency matrix $\btheta_{1:w}^{j*}$ for each motif
$j$ is obtained by drawing its columns $\btheta_k^{j*}$ independently
for $k\in\{1,\ldots,w\}$ from a Dirichlet distribution with parameter
vector $\balpha_1$ (chosen as described below). The background
frequency vector $\btheta_0^*$ is drawn from a Dirichlet distribution
with parameter vector $\balpha_0$. We also define the motif frequency
to be $p_{j} = 0.005$ for each motif $j$ (a typical value in
practice). With these definitions, the data vector $\bS$ is obtained
by drawing each subsequence $\bS_{wi-w+1:wi}$ for $i \in
\{1,\ldots,L/w\}$ from (\ref{EqnGSDef}), using various combinations
of $w$ and $L$. Unlike Assumption~\ref{AssuDeter}, here we use
variable motifs, meaning that $\theta_{k,m}^{j*} \neq1$ for all~$k$
and $m$. We have also done experiments with other values of $p_j$
($p_{j}=0.003$ and $p_{j} = 0.02$), which gave qualitatively the
same results.

We choose $\balpha_0$ and $\balpha_1$ so that the distribution of
$\btheta_0^*$ and that of $\btheta_k^{j*}$ is symmetric in the four
nucleotides;
this means that we must have $\balpha_0 = a_0 \times(1,1,1,1)$ and
$\balpha_1 = a_1 \times(1,1,1,1)$ for some $a_0,a_1 > 0$. Since
motifs are by definition fairly well conserved, we choose $a_1$ so
that the median of $\max_{m \in\{1,\ldots,4\}} \theta_{k,m}^{j*}$
is 0.95
($a_1$ is found numerically). Since background data are
typically more balanced among the four nucleotides, we choose $a_0$ so
that the median of $\max_{m\in\{1,\ldots,4\}} \theta_{0,m}^*$ is 0.3.

For each simulated data vector $\bS$ we run a systematic-scan Gibbs
sampler five times from different initial values and use the
Gelman--Rubin scale factor [\citet{gelmrubi19922}] to detect whether
the chains have converged to different parts of the parameter space,
corresponding to different local modes of the posterior density.
Since the slow mixing in Theorem~\ref{ThmSlowMix2} is caused by
multimodality of the posterior distribution, this approach should
detect the problem effectively.
If different runs of the Markov chain explore different parts of
the parameter space, the Gelman--Rubin scale factor should be large
(typically much larger than $2$), while if they are drawing from the
same distribution the scale factor should be close to~$1$.

In order to detect the worst-case behavior, we take the initial vector
$\bA= (A_1,\ldots, A_{L/w})$ for the first chain to be the vector of
indicators of whether each subsequence $\bS_{wi-w+1:wi}$ was generated\vadjust{\goodbreak}
from motif one. If applicable we initialize the second chain at the
vector of indicators of whether each subsequence $\bS_{wi-w+1:wi}$ was
generated from motif two. The initial vector $\bA$ in all other cases
is generated randomly according to $A_i \iidsim \operatorname{Bernoulli}(p_0)$.
Although in practice one would not know the true motif locations, we
do this to ensure that we detect even very narrow and hard-to-find
modes corresponding to the true motifs. We run each Gibbs
sampler for a burn-in period of $1000$ updates of the entire vector
$\bA$, and then a sampling period of $10,000$ updates of $\bA$. With
these choices standard convergence diagnostics [cf. Geweke
(\citeyear{gewe1992})] that evaluate the convergence of the chains
individually do not detect a convergence problem.

When we run the Gibbs sampler we specify the inference model motif
frequency $p_0$ as $p_0 = \sum_{j=1}^J p_{j}$ (other choices are
investigated below). We specify the prior hyperparameters as
$\beta_{k,m} = 1$ for $k \in\{0,\ldots,w\}$ and $m\in\{1,\ldots,4\}$;
this is the standard choice.

Having simulated the chains, we calculate the
Gelman--Rubin scale factor for the following parameter summaries:
\begin{eqnarray*}
\hat{\theta}_{k,m}(\bA) &\triangleq&\frac{N(\bA^{(k)})_m + \beta_{k,m}} {
|\bN(\bA^{(k)})| + |\bbeta_{k}|},\qquad k \in\{1, \ldots,
w\}, m\in\{1,\ldots,4\},
\\
\hat{\theta}_{0,m}(\bA) &\triangleq&\frac{N(\bA^{c})_m + \beta_{0,m}} {
|\bN(\bA^{c})| + |\bbeta_{0}|}
\end{eqnarray*}
as well as for $|\bA|$, recalling the notation (\ref{EqnCountVecs})
and (\ref{EqnNotation}). The values $\hat{\theta}_{k,m}(\bA)$ and
$\hat{\theta}_{0,m}(\bA)$ are relevant because they are the posterior
means of $\theta_{k,m}$ and $\theta_{0,m}$ given~$\bA$. The posterior
density estimates of $\hat{\theta}_{2,1}(\bA)$ from two different
Markov chains in the case $J=2$ are shown in
Figure~\ref{FigDensPlot}. The Gelman--Rubin scale factor for these
chains is $10.9$, accurately reflecting the fact that the two chains
have converged to different parts of the parameter space.

The top display in Table~\ref{TaboneMotif} addresses the case
of one true motif ($J=1$) and various combinations of $w$ and $L$.
For each combination it reports the percentage out of 20 simulated
datasets for which the maximum Gelman--Rubin scale factor (over the
different parameter summaries) is greater than 1.5. The bottom
display in Table~\ref{TaboneMotif} reports the same quantities for
the case of two true motifs ($J=2$). For one motif no convergence problem
is detected for any of the simulated datasets. For two motifs,
regardless of the value of $w$, there is a severe convergence problem
for large values of $L$.

Finally, we investigate the effect of other choices for $p_0$.
Specifying $p_0 = 0.005$, $p_0= 0.002$ or $p_0 = 0.02$ yields results
that are qualitatively the same as those in Table~\ref{TaboneMotif}.

%t1 #&#
\begin{table}
\def\arraystrech{0.9}
\tablewidth=280pt
\caption{For the cases of one true motif (top) and
two true motifs (bottom), the
percentage of simulated datasets for which the Gelman--Rubin
scale factor from five Gibbs sampling chains is $>1.5$}\label{TaboneMotif}
\begin{tabular*}{280pt}{@{\extracolsep{\fill}}lccc@{}}
\hline
& $\bolds{w=6}$ & $\bolds{w=10}$ & $\bolds{w=15}$ \\
\hline
$J=1$&&&\\
$L/w=2000$ & \phantom{0}0& \phantom{0}0 & \phantom{00}0 \\
$L/w=3000$ & \phantom{0}0 & \phantom{0}0 & \phantom{00}0 \\
$L/w=4000$ &\phantom{0}0& \phantom{0}0 & \phantom{00}0 \\
$L/w=8000$ &\phantom{0}0& \phantom{0}0 & \phantom{00}0 \\[3pt]
$J=2$ & &  & \\
$L/w=2000$ & \phantom{0}0 & 20& \phantom{0}70\\
$L/w=3000$ & 10 & 70&100\\
$L/w=4000$ & 20& 80 &100\\
$L/w=8000$ & 80& 90 &100\\
\hline
\end{tabular*}
\end{table}
\section{\texorpdfstring{Proof of Theorem \protect\ref{ThmSlowMix}}{Proof of Theorem
3.1}}\label{SecModeThm}\vspace*{-3pt}

%s5.1 #&#
\subsection{Specification of the Gibbs sampler}\label{SecGibbs}

Here we give the details of the Gibbs sampler
$T$. Recalling the notation
of Section~\ref{SecMotifFinding},\vadjust{\goodbreak} the sampler iteratively
updates each $A_i$ according to its conditional posterior
distribution, given as follows where $\bA_{[-i]}$ refers to the vector
$\bA$ excluding the $i$th element, where $\bA_{[i,0]}$ is the
vector $\bA$ with the $i$th element replaced by 0, and where
$\bA_{[i,1]}$ is $\bA$ with the $i$th element
replaced by 1. Using (\ref{EqnPostDistn}),
%
%e5.1 #&#
\begin{eqnarray}\label{EqnTrueUpdate}
&&\frac{\pi(A_i=1| \bA_{[-i]}, \bS)} {\pi(A_i=0| \bA_{[-i]}, \bS)}
\nonumber
\\
&&\qquad= \frac{p_0}{1-p_0} \biggl(\frac{\Gamma(\bN(\bA^{c}_{[i,1]}) + \bbeta_0)}{\Gamma(\bN
(\bA^{c}_{[i,0]}) +
\bbeta_0)} \biggr) \frac{\Gamma(|\bN(\bA^{c}_{[i,0]})| + |\bbeta_0|)}{\Gamma(|\bN
(\bA^{c}_{[i,1]})| +
|\bbeta_0|)}
\nonumber\\
&&\qquad\quad{}\times \prod_{k=1}^w \frac{N(\bA^{(k)}_{[i,0]})_{S_{wi-w+k}} +
\beta_{k,S_{wi-w+k}}} {
|\bN(\bA^{(k)}_{[i,0]})| + |\bbeta_{k}|}
\\
&&\qquad= \frac{p_0}{1-p_0} \biggl(\frac{\Gamma(\bN(\bA^{c}_{[i,1]}) + \bbeta_0)}{\Gamma(\bN
(\bA^{c}_{[i,0]}) +
\bbeta_0)} \biggr) \frac{\Gamma(|\bN(\bA^{c}_{[i,0]})| + |\bbeta_0|)}{\Gamma(|\bN
(\bA^{c}_{[i,1]})| +
|\bbeta_0|)} \prod
_{k=1}^w \check{\theta}_{k,S_{wi-w+k}}, \nonumber\\
\eqntext{ i\in\{1,\ldots
L/w\},}
\end{eqnarray}
where the elements of the vector $\check{\btheta}_k$ are the current
estimates of
the frequency of each nucleotide in position $k$ of the motif, that is,
%
%e5.2 #&#
\begin{equation}
\label{EqnThetaK} \check{\theta}_{k,m} \triangleq\frac{N(\bA^{(k)}_{[i,0]})_m +
\beta_{k,m}} {
|\bN(\bA^{(k)}_{[i,0]})| + |\bbeta_{k}|} , \qquad k
\in\{1,\ldots, w\}, m\in\{1,\ldots, M\}.
\end{equation}
For details see \citet{liuneuwlawr1995} and
\citet{jensliuzhou2004}.

%s5.2 #&#
\subsection{\texorpdfstring{Outline of proof of Theorem~\protect\ref{ThmSlowMix}}{Outline of proof of Theorem 3.1}}\label{SecProof31}

Informally, the proof of Theorem~\ref{ThmSlowMix} proceeds by showing
the following results.
\begin{longlist}[\textit{Step} 1.]
\item[\textit{Step} 1.] The order of the spectral gap of the Gibbs sampler is
determined by the unimodality or multimodality of the marginal
posterior distribution of a particular summary vector $\bC(\bA)$ of
$\bA$, denoted by $\bar{\pi}(\bC(\bA) | \bS)$. If
$\bar{\pi}(\bC(\bA)|\bS)$ is multimodal the order of the spectral gap
is determined by the heights of the modes relative to the heights of
the valleys between the modes.

\item[\textit{Step} 2.] When $E \log f(\mathbf{s}| \btheta_{0:w})$ is
multimodal, the marginal posterior distribution $\pi(
\btheta_{0:w} | \bS)$ of the continuous parameters
$\btheta_{0:w}$ is also multimodal, with height of the modes increasing
exponentially in $L$, relative to the height of the valleys between
the modes.

\item[\textit{Step} 3.] The result of step 2 can be mapped to $\bar{\pi}(\bC
(\bA)
| \bS)$, showing that the posterior distribution of $\bC(\bA)$ has
multiple modes with height that grows exponentially in $L$ (relative
to the valleys in between).
\end{longlist}
For simplicity of notation we consider the case $M=2$ (two
nucleotides), although the proof is analogous for any fixed $M$.

Formally, for $\mathbf{s}\in\{1,2\}^w$ any length-$w$ vector of
nucleotides define
%
%e5.3 #&#
\begin{equation}
\label{EqnCADef} C(\bA)_\mathbf{s}\triangleq \bigl\llvert \{i:
A_i=1, \bS_{wi-w+1:wi} = \mathbf{s}\} \bigr\rrvert
\end{equation}
to be the number of instances of motif $\mathbf{s}$ (where $\llvert \{
\ldots\}\rrvert $
indicates the cardinality of a set). Similarly, let
%
%e5.4 #&#
\begin{equation}
\label{EqnCSDef} C(\bS)_\mathbf{s}\triangleq \bigl\llvert \{i\dvtx
\bS_{wi-w+1:wi} = \mathbf {s}\} \bigr\rrvert
\end{equation}
be the number of times that the sequence of nucleotides $\mathbf{s}$
occurs in
the data. Then we must have $C(\bA)_\mathbf{s}\leq C(\bS)_\mathbf
{s}$ for each $\mathbf{s}$,
that is, $\bC(\bA)$ lies in the space
%
%e5.5 #&#
\begin{equation}
\label{EqnCollapseSpace} \bar{\cX} \triangleq\prod_{\mathbf{s}\in\{1,2\}^w}
\bigl\{0, \ldots, C(\bS)_\mathbf{s} \bigr\}.
\end{equation}

The posterior distribution $\pi(\bA| \bS)$ only depends on $\bA$
through $\bC(\bA)$, which can be seen as follows. Using (\ref
{EqnPostDistn}), $\pi(\bA|
\bS)$ depends on $\bA$ via the quantities $|\bA|$, $\bN(\bA^{(k)})$
and $\bN(\bA^c)$. These in turn only depend on $\bC(\bA)$, since
[using (\ref{EqnCountVecs})--(\ref{EqnNotation}) and (\ref{EqnCADef})]
%
%e5.6 #&#
\begin{eqnarray}\label{EqnCountVecsC}
|\bA| & =& \bigl|\bC(\bA)\bigr|,
\nonumber\\
N \bigl(\bA^{(k)} \bigr)_m &=& \sum
_{\mathbf{s}} C(\bA)_\mathbf{s}\indic_{\{
s_k = m\}},\qquad k \in
\{1,\ldots,w\}, m\in\{1,2\},
\\
N \bigl(\bA^c \bigr)_m &=& N(\bS)_m - \sum
_{k=1}^w N \bigl(\bA^{(k)}
\bigr)_m.
\nonumber
\end{eqnarray}

\begin{figure}

\includegraphics{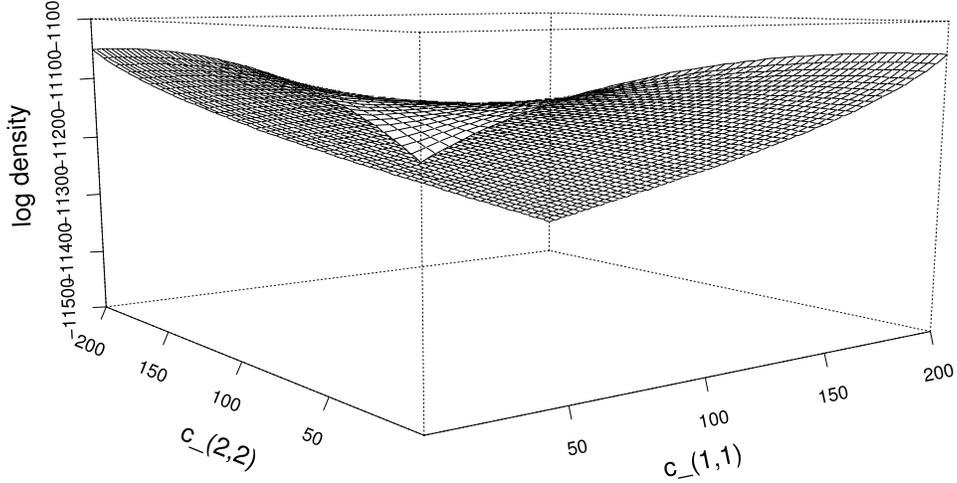}

\caption{The log-density $\log\bar{\pi}$ as a function of
$c_{(1,1)}$ and
$c_{(2,2)}$, fixing $c_{(1,2)}=c_{(2,1)}=0$ and for the case
$w=2$.}\label{FigBimode}
\end{figure}

The marginal posterior distribution of $\bC(\bA)$ is denoted by
%
%e5.7 #&#
\begin{equation}
\label{EqnBarPi} \bar{\pi}(\mathbf{c}|\bS) \triangleq\sum
_{\bA: \bC(\bA) =
\mathbf{c}} \pi( \bA| \bS), \qquad \mathbf{c}\in\bar{\cX}.
\end{equation}
Figure~\ref{FigBimode} illustrates multimodality of $\bar{\pi
}(\mathbf{c}|
\bS)$ for the case $w=2$. For $w=2$ the arguments to the function
$\bar{\pi}$ are the quantities $c_{(1,1)}$, $c_{(1,2)}$, $c_{(2,1)}$
and $c_{(2,2)}$. The data $\bS$ used to create
Figure~\ref{FigBimode} were generated with two true motifs, yielding
the two visible modes of $\bar{\pi}$.

%f3 #&#

We will use Theorem~\ref{ThmDecomp} to bound $\Gap(T)$.
Partition the state space $\cX$ of $T$ according to the value of $\bC
(\bA)$,
%
%e5.8 #&#
\begin{equation}
\label{EqnDCDef} D_{\mathbf{c}} \triangleq \bigl\{\bA\in\cX: \bC(\bA) =
\mathbf{c} \bigr\}, \qquad \mathbf{c}\in\bar{\cX}.
\end{equation}
Define the projection matrix (Appendix~\ref{SecTools}) for $T$ with
respect to
this partition:
%
%e5.9 #&#
\begin{eqnarray}\label{EqnBarT}
\bar{T}(\mathbf{c}_1,\mathbf{c}_2) &\triangleq&
\frac{\sum_{\bA:
\bC(\bA) = \mathbf{c}_1}
\pi(\bA| \bS) T(\bA,D_{\mathbf{c}_{2}})} {
\sum_{\bA: \bC(\bA) = \mathbf{c}_1} \pi(\bA| \bS) }
\nonumber
\\[-8pt]
\\[-8pt]
\nonumber
&=& \sum_{\bA: \bC(\bA) = \mathbf{c}_1} \frac{1}{|D_{\mathbf
{c}_1}|} T(\bA,
D_{\mathbf{c}_2}),\qquad \mathbf {c}_1, \mathbf{c}_2 \in\bar{
\cX}
\end{eqnarray}
since $\pi(\bA| \bS)$ depends on $\bA$ only via $\bC(\bA)$, so that
$\pi(\bA| \bS)$ is equal for all $\bA\in D_{\mathbf{c}_1}$. The
matrix $\bar{T}$ is reversible with respect to $\bar{\pi}$
(Appendix~\ref{SecTools}).

It is easier to obtain useful bounds on $\Gap(T)$
indirectly by relating $T$ to $\bar{T}$ and bounding $\Gap(\bar{T})$
than it is to obtain useful bounds on $\Gap(T)$ directly. The same
technique is utilized in \citet{madrzhen2003} and
\citet{woodschmhube09a}. The
cardinality of the state space $\cX= \{0,1\}^{L/w}$ of $T$ grows
exponentially in $L$ for fixed $w$, while the cardinality of the state
space $\bar{\cX}$ of $\bar{T}$ grows only polynomially in $L$, since [using
(\ref{EqnCollapseSpace})]
%
%e5.10 #&#
\begin{equation}
\label{EqnBarXSize} |\bar{\cX}| \leq(L/w+1)^{2^w}.
\end{equation}
We obtain an upper bound on $\Gap(\bar{T})$ using conductance
(Theorem~\ref{ThmConduct}) and a lower bound on $\Gap(\bar{T})$
using path bounds (Theorem~\ref{ThmPathBound}). Bounds obtained
using these tools can easily be inaccurate by a factor equal to the
cardinality of the space. If we were to obtain these bounds directly
for $\Gap(T)$ they would be loose by an exponential factor in $L$, and
thus unusable for our purposes.

%s5.3 #&#
\subsection{\texorpdfstring{Step 1 of proof of Theorem~\protect\ref{ThmSlowMix}}{Step 1 of proof of Theorem 3.1}}

Consider the graph associated with the reversible matrix $\bar{T}$, with
vertices corresponding to $\mathbf{c}\in\bar{\cX}$ and edges
corresponding to pairs $\mathbf{c}_1, \mathbf{c}_2 \in\bar{\cX}$ having
$\bar{T}(\mathbf{c}_1, \mathbf{c}_2) > 0$. For any $\mathbf
{c}_1,\mathbf{c}_2 \in\bar{\cX}$ let
$\Gamma_{\mathbf{c}_1,\mathbf{c}_2}$ denote the set of paths between
$\mathbf{c}_1$ and
$\mathbf{c}_2$ in the graph that do not have repeated vertices.
Also let $\mathbf{c}\in\gamma$ indicate that the state $\mathbf
{c}\in\bar{\cX}$
is a vertex in the path $\gamma$. Then Theorem~\ref{ThmModes}
formalizes step 1 from Section~\ref{SecProof31}.
%
%th5.1 #&#
\begin{theorem}\label{ThmModes}
$\Gap(T)$ decreases exponentially in $L$ if and only if
%
%e5.11 #&#
\begin{equation}
\label{EqnPathMin} d \triangleq\min_{\mathbf{c}_1,\mathbf{c}_2 \in
\bar{\cX}} \max_{\gamma\in\Gamma_{\mathbf{c}_1,\mathbf{c}_2}}
\min_{\mathbf{c}\in\gamma} \frac{\bar{\pi}(\mathbf{c}|\bS
)}{\bar{\pi}(\mathbf{c}_1|\bS)
\bar{\pi}(\mathbf{c}_2|\bS)}
\end{equation}
decreases exponentially in $L$.
\end{theorem}
The quantity $d$ measures the multimodality of $\bar{\pi}$.
Roughly, think of $\mathbf{c}_1,\mathbf{c}_2$ as being local modes of
$\bar{\pi}$; if all
paths between $\mathbf{c}_1$ and $\mathbf{c}_2$ contain a state with
low probability,
then $d$ is small.

\begin{pf*}{Proof of Theorem~\ref{ThmModes}}
The transition matrix $T$ is nonnegative definite and reversible with
respect to $\pi(\bA|\bS)$. Notice that $T^2$ is also reversible
w.r.t. $\pi(\bA| \bS)$. Using (\ref{EqnDCDef}), let
$T^2|_{D_\mathbf{c}}$
be the restriction of $T^2$ to $D_\mathbf{c}$ as defined in
Appendix~\ref{SecTools}. Then by Lemma~\ref{ThmPowers} and
Theorem~\ref{ThmDecomp},
%
%e5.12 #&#
\begin{eqnarray}\label{EqnDecomp}
\Gap(T)& \geq&\frac{1}{3}\Gap \bigl(T^3 \bigr)  =
\frac{1}{3}\Gap \bigl(T^{1/2} T^2 T^{1/2}
\bigr)
\nonumber
\\[-8pt]
\\[-8pt]
\nonumber
&\geq&\frac{1}{3} \Gap(\bar{T}) \min_{\mathbf{c}\in\bar{\cX} } \Gap
\bigl(T^2|_{D_\mathbf{c}} \bigr).
\end{eqnarray}
Combining with Proposition~\ref{ThmRestBound} below, we
have that
%
%e5.13 #&#
\begin{equation}
\label{EqnGapTInv} \Gap(T)^{-1} = \mathcal{O} \bigl( \Gap(
\bar{T})^{-1} \times L^{2w+3} \bigr).
\end{equation}
Theorem~\ref{ThmDecomp} also gives the bound
$\Gap(T) \leq\Gap(\bar{T}).$
So $\Gap(T)$ is within a polynomial (in $L$) factor of
$\Gap(\bar{T})$.
Theorem~\ref{ThmModes} then follows from
Proposition~\ref{Proppoly2} below.\vadjust{\goodbreak}
\end{pf*}

Finally we give several results used in the proof of Theorem~\ref{ThmModes}.
%
%pr5.1 #&#
\begin{prop}\label{ThmRestBound}
We have
$ [\min_{\mathbf{c}\in\bar{\cX}} \Gap(T^2|_{D_{\mathbf
{c}}}) ]^{-1} = \mathcal{O}(L^{2w+3})$.
\end{prop}

\begin{pf}
Take any $\mathbf{c}\in\bar{\cX}$. For $\mathbf{s}\in\{1,2\}^w$
such that
$C(\bS)_{\mathbf{s}} > 0$ let $\cX_\mathbf{s}\triangleq\{\bz\in\{
0,1\}^{C(\bS)_{\mathbf{s}}}\dvtx \sum_i z_i = c_{\mathbf{s}} \}$. Using (\ref{EqnCSDef}) and
(\ref{EqnDCDef}) the subvector of $\bA\in D_{\mathbf{c}}$ defined by $(A_i\dvtx \bS_{wi-w+1:wi} =\mathbf{s})$
takes values in the space $\cX_\mathbf{s}$.
So there is some bijective map~$h$ such that
%
%e5.14 #&#
\begin{equation}
\label{EqnHMap} \bigl\{h(\bA)\dvtx \bA\in D_{\mathbf{c}} \bigr\} = \prod
_{\mathbf{s}\in
\{1,2\}^w\dvtx C(\bS)_\mathbf{s}> 0} \cX_{\mathbf{s}}.
\end{equation}
Define a transition matrix $\tilde{T}$ having state space $\prod_{\mathbf{s}\in
\{1,2\}^w\dvtx C(\bS)_\mathbf{s}> 0} \cX_{\mathbf{s}}$ and elements
$\tilde{T}(h(\bA),
h(\bA'))\triangleq T^2|_{D_\mathbf{c}}(\bA, \bA')$ for all $\bA,
\bA' \in
D_\mathbf{c}$. Then we have
%
%e5.15 #&#
\begin{equation}
\label{EqnEqualGaps} \Gap(\tilde{T}) = \Gap \bigl(T^2|_{D_\mathbf{c}}
\bigr).
\end{equation}

Using (\ref{EqnDistRest}) and the fact that $T^2$ is reversible
w.r.t. $\pi(\bA| \bS)$, $T^2|_{D_{\mathbf{c}}}$ is reversible
w.r.t. $\pi|_{D_{\mathbf{c}}}(\bA| \bS)$. Since $\pi(\bA| \bS)$
is equal for all
$\bA\in D_{\mathbf{c}}$, $\pi|_{D_{\mathbf{c}}}(\bA| \bS)$ is
uniform on~$D_\mathbf{c}$,
and thus $\tilde{T}$ is also reversible w.r.t. the uniform
distribution.

%The spectral
%gap of $T^2
%|_{D_c}$ is equal to the spectral gap of the Markov transition matrix
%for $h(\bA)$.

We will compare $\tilde{T}$ to a product chain as defined in
Theorem~\ref{ThmProdChain}, denoted by $T^*$. Define $T^{*}$ to have
component chains indexed by $\mathbf{s}\in\{1,2\}^{w}\dvtx C(\bS
)_{\mathbf{s}} > 0$.
The component chains are denoted by $T_\mathbf{s}^{*}$, have state space
$\cX_\mathbf{s}$ and are combined using weights $b_\mathbf{s}=
\frac{C(\bS)_\mathbf{s}}{L/w}$ to form $T^*$. Define $T^{*}_\mathbf
{s}$ to be the
exclusion process with $c_\mathbf{s}$ particles on the complete graph of
$\{1, \ldots, C(\bS)_\mathbf{s}\}$ [\citet{diacsalo1993}]. This
chain is
defined informally as follows, where $\bz\in\cX_\mathbf{s}$ is the current
state. If $c_\mathbf{s}= 0$ or $c_\mathbf{s}= C(\bS)_\mathbf{s}$,
then $\cX_\mathbf{s}$ has a
single state, and the transition matrix $T^{*}_\mathbf{s}$ is trivially
defined. Otherwise, $T_\mathbf{s}^{*}$ chooses an index $j$ uniformly at
random from $\{i \in\{1,\ldots, C(\bS)_{\mathbf{s}}\}\dvtx z_{i} =1\}$.
Then it
chooses an index $\ell\in\{1,\ldots, C(\bS)_{\mathbf{s}}\}$
uniformly at
random. If $z_\ell= 0$, then $z_j$ is changed to $0$, and $z_\ell$ is
changed to $1$; otherwise, the state $\bz$ does not change. The
matrix $T^{*}_\mathbf{s}$ is reversible with respect to the distribution
$\mu_\mathbf{s}$ that is uniform on $\cX_\mathbf{s}$ [\citet
{diacsalo1993}]. So by
Theorem~\ref{ThmProdChain}, $T^{*}$ is reversible with respect to the
uniform distribution on $\prod_{\mathbf{s}\in\{1,2\}^w\dvtx C(\bS
)_\mathbf{s}> 0}
\cX_{\mathbf{s}}$.

By Theorem 3.1 of \citet{diacsalo1993}, for $c_\mathbf{s}>0$ we have
$\Gap(T_\mathbf{s}^{*})\geq1/{c_\mathbf{s}}$, while for $c_\mathbf{s}=0$,
$\Gap(T_\mathbf{s}^{*})=1$. Then Theorem~\ref{ThmProdChain} together with
$c_\mathbf{s}\leq C(\bS)_\mathbf{s}$ yields
%
%e5.16 #&#
\begin{eqnarray}\label{EqnTStar}
\Gap \bigl(T^{*} \bigr) &=& \min_{\mathbf{s}\in\{1,2\}^w\dvtx C(\bS)_{\mathbf{s}}
> 0} b_{\mathbf{s}} \Gap
\bigl(T_\mathbf{s}^* \bigr)
\nonumber
\\[-8pt]
\\[-8pt]
\nonumber
&\geq& \min_{\mathbf{s}\in
\{1,2\}^w\dvtx C(\bS)_{\mathbf{s}}
> 0} \biggl\{ \frac{C(\bS)_\mathbf{s}}{L/w} \biggl( 1 \wedge
\frac{1}{c_\mathbf{s}} \biggr) \biggr\} \geq \frac{w}{L}.
\end{eqnarray}

By (\ref{EqnTransMat}) and the fact that $\bC(\bA) = \bC(\bA') =
\mathbf{c}$
for $\bA,\bA'\in D_\mathbf{c}$, $ T^2(\bA,\bA') > 0$ if and only
if $\exists\mathbf{s}\in
\{1,2\}^w\dvtx C(\bS)_\mathbf{s}> 0$ such that $\bA'$ differs from $\bA
$ only
%are identical except for the following. There are some
%$j,\ell\in\{1,\ldots, L/w\}$ such that $A_j=1$, $A_j'=0$,
%$A_\ell=0$, $A_\ell'=1$ and $\bS_{wj-w+1:wj} = \bS_{w\ell-w+1:w\ell}$.
%So $ T^2|_{D_\bc}(\bA,\bA') > 0$ iff $\exists\bs\in\{1,2\}^w:
%C(\bS)_\bs> 0$ such that the transition from $\bA$ to $\bA'$ only
%changes $\bA$
by swapping two elements of the subvector $(A_{i}\dvtx \bS_{wi-w+1:wi} =
\mathbf{s})$.\vadjust{\goodbreak}
Swapping two elements of a vector in
$\cX_\mathbf{s}$ is precisely the move made by the transition matrix
$T^*_{\mathbf{s}}$. So
%
%e5.17 #&#
\begin{equation}
\label{EqnIFF}  T^{*} \bigl(h(\bA), h \bigl(\bA' \bigr)
\bigr) > 0 \qquad\mbox{iff } T^2 \bigl(\bA, \bA' \bigr)>0
\mbox{ for } \bA, \bA' \in D_\mathbf{c}\dvtx \bA\neq
\bA'.\hspace*{-35pt}
\end{equation}

Using (\ref{EqnTransRest}),
%
%e5.18 #&#
\begin{equation}
\label{EqnT2Equals} T^2|_{D_\mathbf{c}} \bigl(\bA,
\bA' \bigr) = T^2 \bigl(\bA,\bA' \bigr)\qquad
\forall\bA,\bA' \in D_{\mathbf{c}}\dvtx \bA\neq
\bA'.
\end{equation}
Also, by Lemma~\ref{ThmTransProb} below,
%
%e5.19 #&#
\begin{equation}
\label{EqnD1Def} d_1 \triangleq \Bigl[\min_{\bA,\bA'\in\cX: T^2(\bA,\bA') >0}
T^2 \bigl( \bA,\bA' \bigr) \Bigr]^{-1} =
\mathcal{O} \bigl(L^{2w+2} \bigr).
\end{equation}
Combining with (\ref{EqnIFF}) and (\ref{EqnT2Equals}), for any $\bA,\bA'\in D_{\mathbf{c}}$
such that $\bA\neq\bA'$ and $T^{*}(h(\bA),h(\bA')) > 0$,
\[
\tilde{T} \bigl(h(\bA),h \bigl(\bA' \bigr) \bigr) =
T^2 \bigl(\bA,\bA' \bigr) \geq d_1^{-1}
\geq d_1^{-1} T^{*} \bigl(h(\bA),h \bigl(
\bA' \bigr) \bigr).
\]
For $\bA,\bA'\in D_{\mathbf{c}}$ such that $\bA\neq\bA'$ and
$T^{*}(h(\bA),h(\bA')) = 0$ the same inequality holds trivially. So
\[
\tilde{T} \bigl(h(\bA),h \bigl(\bA' \bigr) \bigr) \geq
d_1^{-1} T^{*} \bigl(h(\bA),h \bigl(
\bA' \bigr) \bigr)\qquad\!\!\! \forall h(\bA) \neq h \bigl(\bA'
\bigr) \in\!\prod_{\mathbf{s}\in\{1,2\}^w\dvtx C(\bS)_\mathbf{s}>
0} \!\cX_{\mathbf{s}}.
\]
By Lemma~\ref{moreTransitionLemma} and the fact that both
$\tilde{T}$ and $T^{*}$ are reversible with respect to the uniform
distribution on $\prod_{\mathbf{s}\in\{1,2\}^w\dvtx C(\bS)_\mathbf{s}>
0} \cX_{\mathbf{s}}$,
we then have $\Gap(\tilde{T}) \geq d_1^{-1} \Gap(T^*)$. Combining with
(\ref{EqnEqualGaps}) and (\ref{EqnTStar}),
\begin{eqnarray*}
\Gap \bigl(T^2|_{D_{\mathbf{c}}} \bigr)^{-1} &=& \Gap(
\tilde{T})^{-1} \leq d_1 \Gap \bigl(T^{*}
\bigr)^{-1} \\
&\leq&\frac{d_1
L}{w}
\end{eqnarray*}
regardless of the value of $\mathbf{c}$. By (\ref{EqnD1Def}) $d_1$
does not depend on
$\mathbf{c}$, so\break $\max_{\mathbf{c}\in\bar\cX} \Gap
(T^2|_{D_{\mathbf{c}}})^{-1} =
\mathcal{O}(L^{2w+3}).$
\end{pf}

Lemma~\ref{ThmTransProb} was used in the proof of Proposition~\ref
{ThmRestBound}.
%
%le5.1 #&#
\begin{lemma}\label{ThmTransProb}
We have
\[
\Bigl[\min_{\bA, \bA' \in\cX: T(\bA,\bA') > 0} T \bigl(\bA, \bA' \bigr)
\Bigr]^{-1} = \mathcal{O} \bigl(L^{w+1} \bigr).
\]
\end{lemma}
\begin{pf}
Recall the definition of $\check{\theta}_{k,m}$
from (\ref{EqnThetaK}). Using (\ref{EqnTrueUpdate}),
%
%e5.20 #&#
\begin{eqnarray}
\label{EqnTransProbLower} &&\hspace*{-4pt}\pi(A_i=1 | \bA_{[-i]}, \bS)
\nonumber
\\
&&\hspace*{-4pt}\qquad= \Biggl(p_0 \prod_{k=1}^w
\check{\theta}_{k,S_{wi-w+k}}\Biggr)\nonumber\\[-2pt]
&&\hspace*{-4pt}\qquad\quad{}\bigg/ \Biggl(
p_0 \prod_{k=1}^w \check{\theta}_{k,S_{wi-w+k}} +
(1-p_0)\bigl (\Gamma\bigl(\bN\bigl(\bA^{c}_{[i,0]}\bigr) + \bbeta_0\bigr) \Gamma\bigl(\bigl|\bN
\bigl(\bA^{c}_{[i,1]}\bigr)\bigl| + |\bbeta_0|\bigr)\bigr)\nonumber\\[-2pt]
&&\hspace*{-4pt}\hspace*{170pt}{}/
\bigl(\Gamma\bigl(\bN\bigl(\bA^{c}_{[i,1]}\bigr) + \bbeta_0\bigr) \Gamma\bigl(\bigl|\bN\bigl(\bA^{c}_{[i,0]}\bigr)\bigl| +
|\bbeta_0|\bigr)\bigr)\Biggr)
\\[-2pt]
&&\hspace*{-4pt}\qquad\geq\min \biggl\{ \frac{1}{2}, \Biggl(p_0 \prod_{k=1}^w
\check{\theta}_{k,S_{wi-w+k}}\Biggr)\nonumber\\[-2pt]
&&\hspace*{-4pt}\hspace*{71pt}{}\bigg/\bigl(
2(1-p_0)\bigl(\Gamma\bigl(\bN\bigl(\bA^{c}_{[i,0]}\bigr) + \bbeta_0\bigr) \Gamma
\bigl(\bigl|\bN\bigl(\bA^{c}_{[i,1]}\bigr)\bigr| + |\bbeta_0|\bigr)\bigr)\nonumber\\[-2pt]
&&\hspace*{-4pt}\hspace*{124pt}{}/
\bigl(\Gamma\bigl(\bN\bigl(\bA^{c}_{[i,1]}\bigr) + \bbeta_0\bigr) \Gamma\bigl(\bigl|\bN\bigl(\bA^{c}_{[i,0]}\bigr)\bigr| +
|\bbeta_0|\bigr)\bigr)\bigr) \biggr\}.\nonumber
\end{eqnarray}
Also, by (\ref{EqnCountVecs}) and the definitions of $\bA_{[i,0]}$
and $\bA_{[i,1]}$,
\begin{eqnarray*}
N \bigl(\bA^c_{[i,0]} \bigr)_m &=& N \bigl(
\bA^c_{[i,1]} \bigr)_m + \sum
_{k=1}^w \indic_{\{\bS_{wi-w+k} = m\}}, \qquad m\in\{1,2\},
\\[-2pt]
\bigl|\bN \bigl(\bA^c_{[i,0]} \bigr)\bigr| &=&\bigl |\bN \bigl(
\bA^c_{[i,1]} \bigr)\bigr| + w.
\end{eqnarray*}
So
\begin{eqnarray*}
&&\frac{\Gamma(\bN(\bA^{c}_{[i,0]}) + \bbeta_0) \Gamma(|\bN(\bA^{c}_{[i,1]})| + |\bbeta_0|)} {
\Gamma(\bN(\bA^{c}_{[i,1]}) + \bbeta_0) \Gamma(|\bN(\bA^{c}_{[i,0]})| +
|\bbeta_0|)}
\\[-2pt]
&&\qquad=\Biggl(\prod_{m=1}^2  \bigl(N\bigl(\bA^{c}_{[i,1]}\bigr)_m +
\beta_{0,m} \bigr) \bigl(N\bigl(\bA^{c}_{[i,1]}\bigr)_m + \beta_{0,m}+1
\bigr)\\[-2pt]
&&\hspace*{56pt}{} \cdots\bigl(N\bigl(\bA^{c}_{[i,1]}\bigr)_m
+ \beta_{0,m} + \sum_{k=1}^w \indic_{\{S_{wi-w+k}=m\}} -1 \bigr)\Biggr)\\[-2pt]
&&\qquad\quad{}/\bigl(
\bigl(\bigl|\bN\bigl(\bA^{c}_{[i,1]}\bigr)\bigr| +
|\bbeta_{0}|\bigr)\bigl (\bigl|\bN\bigl(\bA^{c}_{[i,1]}\bigr)\bigr| +
|\bbeta_{0}| + 1\bigr) \\[-2pt]
&&\hspace*{90pt}{}\cdots\bigl(\bigl|\bN\bigl(\bA^{c}_{[i,1]}\bigr)\bigr| +
|\bbeta_{0}|+w-1\bigr)\bigr)
\\[-2pt]
&&\qquad\leq\Biggl(\prod_{m=1}^2  \bigl(\bigl|\bN\bigl(\bA^{c}_{[i,1]}\bigr)\bigr| +
|\bbeta_{0}| \bigr) \bigl(\bigl|\bN\bigl(\bA^{c}_{[i,1]}\bigr)\bigr| + |\bbeta_{0}|+1
\bigr)\\[-2pt]
&&\hspace*{56pt}{} \cdots\Biggl(\bigl|\bN\bigl(\bA^{c}_{[i,1]}\bigr)\bigr|
+ |\bbeta_{0}| + \sum_{k=1}^w \indic_{\{S_{wi-w+k}=m\}} -1 \Biggr)\Biggr) \\[-2pt]
&&\qquad\quad{}/
\bigl(\bigl(\bigl|\bN\bigl(\bA^{c}_{[i,1]}\bigr)\bigr| +
|\bbeta_{0}|\bigr) \bigl(\bigl|\bN\bigl(\bA^{c}_{[i,1]}\bigr)\bigr| +
|\bbeta_{0}| + 1\bigr) \\[-2pt]
&&\hspace*{90pt}{}\cdots\bigl(\bigl|\bN\bigl(\bA^{c}_{[i,1]}\bigr)\bigr| +
|\bbeta_{0}|+w-1\bigr)\bigr)
\\
&&\qquad\leq1.
\end{eqnarray*}
%
%There are $w$ factors in the numerator and in the denominator of
%(\ref{EqnRatioLowerBd}), and each factor in the numerator can be
%paired with a factor in the denominator that is at least as large. So
%the ratio
%in (\ref{EqnRatioLowerBd}) is $\leq1$.
Combining with (\ref{EqnThetaK}) and
(\ref{EqnTransProbLower}),
\[
\pi(A_i=1 | \bA_{[-i]}, \bS) \geq\frac{p_0 \prod_{k=1}^w
\check{\theta}_{k,S_{wi-w+k}}}{2} \geq
\frac{p_0}{2} \prod_{k=1}^w
\frac{\beta_{k,S_{wi-w+k}}}{L+|\bbeta_k| },
\]
which does not depend on $\bA$ or $i$. So
$ [ \min_{\bA, i} \pi(A_i=1 | \bA_{[-i]}, \bS)  ]^{-1}
= \mathcal{O}(L^w).
$
Analogously,
$ [ \min_{\bA, i} \pi(A_i=0 | \bA_{[-i]}, \bS)  ]^{-1}
= \mathcal{O}(L^w).$
Using (\ref{EqnTransMat}) then yields the desired result.
\end{pf}

Proposition~\ref{Proppoly2} was used in the proof of Theorem~\ref{ThmModes}.
%
%pr5.2 #&#
\begin{prop} \label{Proppoly2}
$\Gap(\bar{T})$ is within a polynomial (in $L$) factor of $d$.
Specifically,
$\Gap(\bar{T}) = \mathcal{O}(d \times L^{2^w})$ and
$\Gap(\bar{T})^{-1} = \mathcal{O} ( d^{-1} \times
L^{w+1+2^{w+1} + 2^w}  )$.
\end{prop}
The bounds in Proposition~\ref{Proppoly2} only rely on
Lemma~\ref{ThmTransProb} and the fact (\ref{EqnBarXSize}) that
$|\bar{\cX}|$ grows polynomially in $L$, and could be improved by
leveraging additional properties of $\bar{T}$ (at the cost of some
technical complexity). However, Proposition~\ref{Proppoly2} is
sufficient for our purposes.

\begin{pf*}{Proof of Proposition~\ref{Proppoly2}}

\textit{Upper bound}.
Suppress the dependence\vspace*{1pt} of $\bar{\pi}(\mathbf{c}|\bS)$ on $\bS$ for
simplicity of notation. Let $\mathbf{c}_1,\mathbf{c}_2 \in\bar{\cX
}$ be a pair of
states that achieve the minimum in definition (\ref{EqnPathMin})
of $d$, so that $d = \max_{\gamma\in\Gamma_{\mathbf{c}_1,\mathbf{c}_2}}
\min_{\mathbf{c}\in\gamma} \frac{\bar{\pi}(\mathbf{c})}{\bar
{\pi}(\mathbf{c}_1)
\bar{\pi}(\mathbf{c}_2)}.$ For any $\gamma\in\Gamma_{\mathbf
{c}_1,\mathbf{c}_2}$, let
$\mathbf{c}_\gamma\triangleq\operatorname{argmin}_{\mathbf{c}\in\gamma}
\frac{\bar{\pi}(\mathbf{c})}{\bar{\pi}(\mathbf{c}_1)\bar{\pi
}(\mathbf{c}_2)}$ (in the case\vspace*{1pt}
of a tie choose the state earliest in the path). Defining the set $E
= \{\mathbf{c}_\gamma: \gamma\in\Gamma_{\mathbf{c}_1,\mathbf
{c}_2}\}$, we have
%
%e5.21 #&#
\begin{equation}
d = \max_{\mathbf{c}_{\gamma} \in E} \frac{\bar{\pi}(\mathbf
{c}_\gamma)}{\bar{\pi}(\mathbf{c}_1)
\bar{\pi}(\mathbf{c}_2)}.\label{EqnDAlt}
\end{equation}
The set $E$ separates $\mathbf{c}_1$ and $\mathbf{c}_2$ in the sense
that there is
no path $\gamma\in\Gamma_{\mathbf{c}_1,\mathbf{c}_2}$ that does
not include some
state in $E$. If $\mathbf{c}_1 \in E$, then there is some $\gamma\in
\Gamma_{\mathbf{c}_1,\mathbf{c}_2}$ for which $\mathbf{c}_{\gamma}
= \mathbf{c}_1$, and so $d \geq
\frac{1}{\bar{\pi}(\mathbf{c}_2)} \geq1$. In this case $\Gap(\bar
{T}) \leq
2 d (L/w+1)^{2^w}$ holds since $\Gap(\bar{T}) \leq2$.

Now consider the case $\mathbf{c}_1 \notin E$. Let $B$ be the set of states
$\mathbf{c}\in\bar{\cX}$ that are not reachable from $\mathbf
{c}_1$ without going
through $E$,
\[
B \triangleq \bigl\{\mathbf{c}\in\bar{\cX}\dvtx \forall\gamma\in
\Gamma_{\mathbf{c}_1,\mathbf{c}} \mbox{ there is some } \mathbf{c}' \in\gamma
\mbox{ s.t. } \mathbf {c}' \in E \bigr\}.
\]
We have that $\mathbf{c}_2 \in B$, and $\mathbf{c}_1 \in B^c$ since
$\mathbf{c}_1 \notin E$.
Also, the only states $\mathbf{c}\in B$ for which $\bar{T}(\mathbf
{c},B^c)>0$ satisfy
$\mathbf{c}\in E$, which can be seen as follows. Otherwise, $\exists\mathbf
{c}\in B\setminus E$
and $\mathbf{c}_3 \in B^c$ for which $\bar{T}(\mathbf{c},\mathbf
{c}_3)>0$. Since $\mathbf{c}_3 \in B^c$,
there is a path $\gamma\in\Gamma_{\mathbf{c}_1,\mathbf{c}_3}$ that
does not go through
$E$. But since $\bar{T}(\mathbf{c}_3,\mathbf{c}) >0$ ($\bar{T}$ is
reversible), there
is also a path $\gamma\in\Gamma_{\mathbf{c}_1,\mathbf{c}}$ that
does not go through
$E$, which is a contradiction.

Using these facts, (\ref{EqnBarXSize}) and (\ref{EqnDAlt}), the
conductance of $B$ (Theorem~\ref{ThmConduct}) is
\begin{eqnarray*}
\Phi_{\bar{T}}(B) &= &\frac{\sum_{\mathbf{c}\in B} \bar{\pi
}(\mathbf{c})
\bar{T}(\mathbf{c},B^c)}{ \bar{\pi}(B) \bar{\pi}(B^c)} \leq\frac
{\sum_{\mathbf{c}\in E} \bar{\pi}(\mathbf{c})
\bar{T}(\mathbf{c},B^c)}{ \bar{\pi}(B) \bar{\pi}(B^c)}
\\
&\leq&\frac{\sum_{\mathbf{c}\in E} \bar{\pi}(\mathbf{c})}{\bar
{\pi}(B)
\bar{\pi}(B^c)} \leq\frac{\sum_{\mathbf{c}\in E} \bar{\pi
}(\mathbf{c})}{\bar{\pi}(\mathbf{c}_1)
\bar{\pi}(\mathbf{c}_2)}
\\
&\leq&|E| \max_{\mathbf{c}\in E} \frac{\bar{\pi}(\mathbf
{c})}{\bar{\pi}(\mathbf{c}_1)
\bar{\pi}(\mathbf{c}_2)} = |E| d \leq|\bar{\cX}| d \leq
d (L/w+1)^{2^w}.
\end{eqnarray*}
Using Theorem~\ref{ThmConduct}, $\Gap(\bar{T}) \leq2 \Phi_{\bar
{T}}(B) \leq2 d (L/w+1)^{2^w}$ as
claimed.

\textit{Lower bound.}
Recall that the transition matrix $\bar{T}$ has state space
$\bar{\cX}$ and is reversible with respect to $\bar{\pi}$. Using
(\ref{EqnDCDef})--(\ref{EqnBarT}), for $\mathbf{c},\mathbf{c}' \in
\bar{\cX}$ such
that $\sum_\mathbf{s}|c_\mathbf{s}-c'_\mathbf{s}| \leq1$, $\bar
{T}(\mathbf{c},\mathbf{c}') > 0$ and
otherwise $\bar{T}(\mathbf{c},\mathbf{c}')=0$. Also
%
%e5.22 #&#
\begin{equation}
\label{EqnbarT} \bar{T} \bigl(\mathbf{c},\mathbf{c}' \bigr) \geq
\min_{\bA\in D_{\mathbf{c}}} T(\bA,D_{\mathbf{c}'}) \geq\min_{\bA\in D_{\mathbf{c}}}
\max_{\bA' \in D_{\mathbf{c}'}} T \bigl(\bA,\bA' \bigr)\qquad \forall\mathbf{c},
\mathbf{c}'\in\bar{\cX}.
\end{equation}
If
$\bar{T}(\mathbf{c},\mathbf{c}') > 0$, then for every $\bA\in
D_{\mathbf{c}}$ there is
some $\bA' \in D_{\mathbf{c}'}$ with $T(\bA,\bA') > 0$. By
Lemma~\ref{ThmTransProb} and (\ref{EqnbarT}),
%
%e5.23 #&#
\begin{equation}
\label{EqnbarTmin} \Bigl[\min_{\mathbf{c},\mathbf{c}'\in\bar{\cX}\dvtx \bar{T}(\mathbf
{c},\mathbf{c}') > 0} \bar{T} \bigl(\mathbf{c},
\mathbf{c}' \bigr) \Bigr]^{-1}= \mathcal{O}
\bigl(L^{w+1} \bigr).
\end{equation}

We will use Theorem~\ref{ThmPathBound} to obtain a bound for
$\Gap(\bar{T})$; let $\E$ be the set of edges in the graph of
$\bar{T}$. For $(z,v) \in\E$ and $\gamma$ a path in the graph, let
$(z,v) \in\gamma$ indicate that the edge $(z,v)$ is in the path
$\gamma$ (as distinct from $z\in\gamma$ which indicates that the
vertex $z$ is in $\gamma$). To apply Theorem~\ref{ThmPathBound} we
need to define a path $\gamma_{x,y}$ for every pair of states $x,y\in
\bar{\cX}$. Suppressing the dependence of $\bar{\pi}$
on $\bS$, choose $\gamma_{x,y}$ to be any path that maximizes
$\min_{(z,v)\in\gamma} \frac{\bar{\pi}(z)}{\bar{\pi}(x)
\bar{\pi}(y)}$. Letting $|\{\ldots\}|$ denote the cardinality of a
set, the path constant $\rho$ defined in Theorem~\ref{ThmPathBound}
satisfies
\begin{eqnarray*}
\rho&=& \max_{(z,v)\in\E} \frac{1}{\bar{\pi}(z)\bar{T}(z,v)} \sum
_{\gamma_{x,y}\ni(z,v) } \bar{\pi}(x)\bar{\pi}(y) \operatorname{len}(
\gamma_{x,y})
\\
&\leq&\frac{1}{\min_{(z,v)\in\E} \bar{T}(z,v)} \biggl[\max_{(z,v)\in\E} \frac{1}{\bar{\pi}(z)} \sum
_{\gamma
_{x,y}\ni(z,v) } \bar{\pi}(x)\bar{\pi}(y) \biggr]
\max_{x,y} \operatorname{len}(\gamma_{x,y})
\\
&\leq&\frac{1}{\min_{(z,v)\in\E} \bar{T}(z,v)} \Bigl[\max_{(z,v)\in\E} \bigl| \bigl\{
\gamma_{x,y} \ni(z,v) \bigr\}\bigr| \Bigr] \\
&&{}\times\biggl[\max_{(z,v)\in\E}
\max_{\gamma_{x,y} \ni(z,v)} \frac{\bar
{\pi}(x) \bar{\pi}(y)}{\bar{\pi}(z)} \biggr] \max_{x,y} \operatorname{len}(
\gamma_{x,y})
\\
&= &\frac{1}{\min_{(z,v)\in\E} \bar{T}(z,v)} \Bigl[\max_{(z,v)\in\E}\bigl | \bigl\{\gamma_{x,y}
\ni(z,v) \bigr\}\bigr| \Bigr]\\
&&{}\times \biggl[\min_{x,y \in\bar{\cX}} \min_{(z,v) \in\gamma_{x,y}}
\frac{\bar{\pi}(z)}{\bar{\pi}(x) \bar{\pi}(y)} \biggr]^{-1} \max_{x,y} \operatorname{len}(
\gamma_{x,y})
\\
&=& \frac{1}{\min_{(z,v)\in\E} \bar{T}(z,v)} \Bigl[\max_{(z,v)\in\E} \bigl| \bigl\{\gamma_{x,y}
\ni(z,v) \bigr\}\bigr| \Bigr]\\
&&{}\times \biggl[\min_{x,y \in\bar{\cX}} \max_{\gamma\in\Gamma_{x,y}}
\min_{(z,v) \in\gamma} \frac{\bar{\pi}(z)}{\bar{\pi}(x) \bar
{\pi}(y)} \biggr]^{-1}
\max_{x,y} \operatorname{len}(\gamma_{x,y})
\\
&\leq&\frac{1}{\min_{(z,v)\in\E} \bar{T}(z,v)} \Bigl[\max_{(z,v)\in\E}\bigl | \bigl\{
\gamma_{x,y} \ni(z,v) \bigr\}\bigr| \Bigr]\\
&&{}\times \biggl[\min_{x,y \in\bar{\cX}}
\max_{\gamma\in\Gamma_{x,y}} \min_{z \in\gamma} \frac{\bar{\pi}(z)}{\bar{\pi}(x) \bar{\pi
}(y)} \biggr]^{-1}
\max_{x,y} \operatorname{len}(\gamma_{x,y}).
\end{eqnarray*}
From (\ref{EqnBarXSize}) $|\bar{\cX}| = \mathcal{O}(L^{2^w})$, so the
maximum length of paths is\break $ \max_{x,y} \operatorname{len}(\gamma_{x,y}) =
\mathcal{O}(L^{2^w})$, and the total
number of paths is no more than $|\bar{\cX}|^2 =
\mathcal{O}(L^{2^{w+1}})$. By Theorem~\ref{ThmPathBound} and (\ref
{EqnbarTmin}),
$\Gap(\bar{T})^{-1} \leq\rho=\break \mathcal{O} (d^{-1}
L^{w+1+2^{w+1}+2^w}  ).$
\end{pf*}

%s5.4 #&#
\subsection{\texorpdfstring{Step 2 of proof of Theorem~\protect\ref{ThmSlowMix}}{Step 2 of proof of Theorem 3.1}}

Recalling
that $M=2$, there are $w+1$ free parameters $\theta_{k,1} \in[0,1]$
for $k \in\{0,\ldots,w\}$, so we write
%
%e5.24 #&#
\begin{equation}
\btheta_{0:w} \in [0,1]^{w+1}. \label{EqnLambdaDef}
\end{equation}
Theorem~\ref{ThmMultimode} formalizes step 2, using Theorem~\ref{ThmBerk}.

%th5.2 #&#
\begin{theorem}\label{ThmMultimode}
Under Assumption~\ref{AssuGener}, if $E \log f(\mathbf{s}| \btheta_{0:w})$
is multimodal, then there exist $\varepsilon> 0$ and two
sets $B_1, B_2 \subset[0,1]^{w+1}$ separated by Euclidean distance
$\varepsilon$ such that
%
%e5.25 #&#
\begin{equation}
\frac{\pi (\btheta_{0:w} \notin B_1 \cup B_2|
\bS )}{\pi (\btheta_{0:w} \in B_1 | \bS )}\quad \mbox{and}\quad \frac{\pi (\btheta_{0:w} \notin B_1 \cup B_2|
\bS )}{\pi (\btheta_{0:w} \in B_2 | \bS )}\label {EqnSetRatio}
\end{equation}
decrease exponentially in $L$, almost surely.
\end{theorem}
\begin{pf}
The inference model assumes $\bS_{wi-w+1:wi} \iidsim
f(\mathbf{s}|\btheta_{0:w})$ for $i \in\{1,\ldots, L/w\}$. By
Assumption~\ref{AssuGener}, the generative model assumes\break
$\bS_{wi-w+1:wi}\iidsim g(\mathbf{s})$, fitting into the framework of
Theorem~\ref{ThmBerk}. Using the notation of that theorem, consider
the case where $\eta(\btheta_{0:w}) = E \log f(\mathbf{s}| \btheta_{0:w})$
is multimodal. Then there are $\btheta_{0:w}^{(j)} \in[0,1]^{w+1}$
and $F_j \subset[0,1]^{w+1}$ for $j \in\{1,2\}$ such that
$\btheta_{0:w}^{(j)} \in F_j$ and~(\ref{EqnMultimodeDef}) holds. Let
$\phi_j \triangleq\sup_{\partial F_j} \eta(\btheta_{0:w}) <
\eta(\btheta_{0:w}^{(j)})$ for $j \in\{1,2\}$, and take any $\xi$ for
which
%
%e5.26 #&#
\begin{equation}
\xi\in \Bigl( \min_{j\in\{1,2\}} \phi_j, \min_{j \in
\{1,2\}}\eta
\bigl(\btheta_{0:w}^{(j)} \bigr) \Bigr). \label{EqnXiDef}
\end{equation}
So $\xi> \phi_j$ for some $j \in\{1,2\}$; assume WLOG that $\xi
>\phi_1$.\vadjust{\goodbreak}

Define the sets
%
%e5.27 #&#
\begin{eqnarray}\label{EqnTildeV}
V &\triangleq& \bigl\{\btheta_{0:w} \in [0,1]^{w+1}: \eta(
\btheta_{0:w}) \geq\xi \bigr\},
\nonumber
\\[-8pt]
\\[-8pt]
\nonumber
\tilde{B}_1 &\triangleq& F_1 \cap V,\qquad
\tilde{B}_2 \triangleq%V \setminus\tilde{B}_1 =
V \setminus F_1.
\end{eqnarray}
By
(\ref{EqnXiDef}), $\eta(\btheta_{0:w}^{(1)}) > \xi$ and
$\btheta_{0:w}^{(1)} \in\tilde{B}_1$. Also using
(\ref{EqnMultimodeDef}), $\eta(\btheta_{0:w}^{(2)}) > \xi$ and
$\btheta^{(2)}_{0:w} \in F_2 \cap V \subset\tilde{B}_2$, so
%
%e5.28 #&#
\begin{equation}
V = \tilde{B}_1 \cup\tilde{B}_2 \quad\mbox{and}\quad
\sup_{\btheta_{0:w}
\in
\tilde{B}_j} \eta( \btheta_{0:w}) > \xi ,\qquad j \in\{1,2\}.
\label{EqnBTilde}
\end{equation}

If $\nexists\varepsilon> 0$ such that $\tilde{B}_1$ and
$\tilde{B}_2$ are separated by distance $\varepsilon$, then (since
$[0,1]^{w+1}$ is compact) $\operatorname{cl}(\tilde{B}_1) \cap
\operatorname{cl}(\tilde{B}_2) \neq\varnothing$. By (\ref{EqnTildeV})
$\tilde{B}_1 \subset
F_1$ and $\tilde{B}_2 \subset[0,1]^{w+1} \setminus F_1$ so
$\operatorname{cl}(\tilde{B}_1) \cap\operatorname{cl}(\tilde{B}_2) \subset
\partial F_1$. This is a contradiction since, due to
(\ref{EqnTildeV})--(\ref{EqnBTilde}) and the continuity of $\eta$,
\[
\inf_{\operatorname{cl}(\tilde{B}_1) \cap\operatorname{cl}(\tilde{B}_2)} \eta \geq\inf_{\operatorname{cl}(\tilde{B}_1) \cup\operatorname{cl}(\tilde{B}_2)} \eta= \inf_{\tilde{B}_1 \cup\tilde{B}_2}
\eta\geq\xi> \phi_1 = \sup_{\partial
F_1} \eta.
\]
So $\exists\varepsilon> 0$ such that $\tilde{B}_1$ and $\tilde{B}_2$
are separated by distance $\varepsilon$.

In order to satisfy assumption~(2) of Theorem~\ref{ThmBerk} we remove points
from the space $[0,1]^{w+1}$ of $\btheta_{0:w}$ for which
$\exists\mathbf{s}\in\{1,2\}^w\dvtx f(\mathbf{s}| \btheta_{0:w}) = 0$. This
results in the space
%
%e5.29 #&#
\begin{equation}
\Lambda= \bigl((0,1)\times[0,1]^w \bigr) \cup \bigl( [0,1] \times
(0,1)^w \bigr).\label{EqnRealLambda}
\end{equation}
The fact that $f(\mathbf{s}| \btheta_{0:w})>0$ for all $\mathbf{s}$
and all
$\btheta_{0:w} \in\Lambda$ is a consequence of (\ref{EqnFSDef}): if
$\theta_{0,1} \in(0,1)$, then
$f(\mathbf{s}| \btheta_{0:w}) > 0$ for all $\mathbf{s}$, and the
same holds
when $\theta_{k,1} \in(0,1)$ for all $k \in\{1,\ldots, w\}$.
%Conversely, if $\theta_{0,1} \in\{0,1\}$ and $\exists k \in
%some $\bs\in\{1,2\}^w$ such that $f(\bs| \btheta_{0:w}) = 0$.
%For example, if $\theta_{0,1} = 0$ and $\theta_{1,1} = 1$ then
%$f\left( (2,1,\ldots,1) | \btheta_{0:w}\right) = 0$.

Taking
%
%e5.30 #&#
\begin{equation}
B_j \triangleq\tilde{B}_j \cap \Lambda, \qquad j \in\{1,2\}
\label{EqnBDef}
\end{equation}
$B_1$ and $B_2$ are separated by distance $\varepsilon$. Due to
(\ref{EqnMultimodeDef}), $\btheta_{0:w}^{(j)}$ is not a limit point
of $[0,1]^{w+1} \setminus F_j$ for $j \in\{1,2\}$. Using (\ref
{EqnRealLambda}) there
are points $\btheta_{0:w} \in\Lambda$ arbitrarily close to
$\btheta_{0:w}^{(j)} \in[0,1]^{w+1}$. From (\ref{EqnXiDef}) and the
continuity of $\eta$, all such points $\btheta_{0:w}$ close enough to
$\btheta_{0:w}^{(j)}$ have $\btheta_{0:w} \in F_j \cap\Lambda$ and
$\eta(\btheta_{0:w}) > \xi$. Using (\ref{EqnMultimodeDef}),
(\ref{EqnTildeV}) and~(\ref{EqnBDef}) these points are in $B_j$, so
%
%e5.31 #&#
\begin{equation}
\sup_{\btheta_{0:w} \in B_j} \eta( \btheta_{0:w}) > \xi,\qquad  j \in\{1,2\}.
\label{EqnBProps}
\end{equation}

Let $\operatorname{Int}(\cdot)$ denote set interior with respect to the space
$\Lambda$. We claim that
%
%e5.32 #&#
\begin{equation}
\bigl\{\btheta_{0:w} \in\Lambda: \eta(\btheta_{0:w}) > \xi
\bigr\} \subset \operatorname{Int}(B_1) \cup\operatorname{Int}(B_2),
\label{EqnInBs}
\end{equation}
which can be seen as follows. Take any $\btheta_{0:w} \in\Lambda$
such that $\eta(\btheta_{0:w}) > \xi$. By
(\ref{EqnTildeV}), (\ref{EqnBTilde}), (\ref{EqnBDef}) and since
$\eta$ is continuous, $\btheta_{0:w} \in\operatorname{Int}(B_1 \cup B_2)$.
Because $B_1$ and $B_2$ are separated by distance $\varepsilon>0$,
$\btheta_{0:w} \in\operatorname{Int}(B_1) \cup\operatorname{Int}(B_2)$.

Define the alternative parameter spaces $\Lambda_1 \triangleq\Lambda
\setminus\operatorname{Int}(B_2)$ and $\Lambda_2 \triangleq\Lambda
\setminus
\operatorname{Int}( B_1)$. By (\ref{EqnBProps}), $\sup_{\Lambda_j} \eta>
\xi$ for $j \in\{1,2\}$. So $\delta_j \triangleq
\frac{\sup_{\Lambda_j}\eta- \xi}{2} > 0$ for $j \in\{1,2\}$.
Combining (\ref{EqnInBs}) with the fact that $\sup_{\Lambda_j} \eta-
\delta_j > \xi$,
%
%e5.33 #&#
\begin{eqnarray}
\label{EqnUDef} U^j_{\delta_j} \triangleq \Bigl\{
\btheta_{0:w} \in \Lambda_j \dvtx \eta(\btheta_{0:w})
\geq \sup_{\Lambda_j} \eta - \delta_j \Bigr\} \subset
\operatorname{Int}(B_j) \subset B_j,
\nonumber
\\[-8pt]
\\[-8pt]
  \eqntext{j \in\{1,2\}.}
\end{eqnarray}

By (\ref{EqnUDef}),
%
%e5.34 #&#
\begin{equation}
\qquad\Lambda\setminus(B_1 \cup B_2) \subset \Lambda
\setminus \bigl(B_1 \cup\operatorname{Int}(B_2) \bigr)
\label{EqnIsSubset} \subset \Lambda\setminus \bigl(U^1_{\delta_1}
\cup\operatorname{Int}(B_2) \bigr) = \Lambda_1 \setminus
U^1_{\delta_1}.
\end{equation}
Analogously, $\Lambda\setminus(B_1 \cup B_2) \subset\Lambda_2
\setminus U^2_{\delta_2}$.

The regularity conditions of Theorem~\ref{ThmBerk} are verified in
the supplementary material [\citet{woodrose12}] for each of the
parameter spaces $\Lambda_1$ and~$\Lambda_2$. We apply that theorem for each of $j \in\{1,2\}$, with
parameter space $\Lambda_j$,
using $\delta= \delta_j$ and taking $n=L/w$. This yields
\[
\limsup_{n \rightarrow\infty} \biggl( \frac{P_{n}(\Lambda_j\setminus U^j_{\delta_j} )}{P_{n}(U^j_{\delta_j})} \biggr)^{1/n} \leq
e^{-\delta_j} \qquad\mbox{a.s. } j \in\{1,2\}.
\]
Combining with (\ref{EqnUDef})--(\ref{EqnIsSubset}) and the fact that
$[0,1]^{w+1}\setminus\Lambda$ has probability zero under
$\pi(\btheta_{0:w} | \bS)$, for $j \in\{1,2\}$
\begin{eqnarray*}
&&\limsup_{L \rightarrow\infty} \biggl( \frac{\pi (\btheta_{0:w} \in [0,1]^{w+1} \setminus(B_1 \cup B_2)
|   \bS )} {
\pi(\btheta_{0:w} \in B_j | \bS)} \biggr)^{{1}/{(L/w)}}
\\
&&\qquad= \limsup_{L \rightarrow\infty} \biggl( \frac{\pi (\btheta_{0:w} \in\Lambda\setminus(B_1 \cup B_2)
|  \bS )} {
\pi(\btheta_{0:w} \in B_j | \bS)} \biggr)^{{1}/{(L/w)}}
\\
&&\qquad\leq \limsup_{L \rightarrow\infty} \biggl( \frac{\pi(\btheta_{0:w} \in\Lambda_j \setminus U_{\delta_j}^j |
\bS)} {
\pi(\btheta_{0:w} \in U_{\delta_j}^j | \bS)} \biggr)^{{1}/{(L/w)}}\\
&&\qquad=
\limsup_{n \rightarrow\infty} \biggl( \frac{P_{n}(\Lambda_j\setminus U^j_{\delta_j} )}{P_{n}(U^j_{\delta_j})} \biggr)^{1/n} \leq
e^{-\delta_j}
\end{eqnarray*}
almost surely.
\end{pf}

%s5.5 #&#
\subsection{\texorpdfstring{Step 3 of proof of Theorem~\protect\ref{ThmSlowMix}}{Step 3 of proof of Theorem 3.1}}

Finally we formalize step 3.

%th5.3 #&#
\begin{theorem}\label{ThmMapSpaces}
If there exist $\varepsilon> 0$ and two
sets $B_1, B_2 \subset[0,1]^{w+1}$ separated by Euclidean distance
$\varepsilon$ such that the ratios in (\ref{EqnSetRatio})
decrease exponentially in $L$, then the quantity $d$ in (\ref{EqnPathMin})
decreases exponentially in~$L$.
\end{theorem}
Theorem~\ref{ThmMapSpaces} is proven in the supplementary material [\citet{woodrose12}].
Theorems~\ref{ThmModes}--\ref{ThmMapSpaces} together imply
Theorem~\ref{ThmSlowMix}.

%s6 #&#
\section{Conclusions}\label{SecConclude}

The Gibbs sampling method is a popular approach to finding gene
regulatory binding motifs, but its poor convergence in practice means
that it can only be used to generate candidate motifs that must be
ranked using a secondary criterion. If one could efficiently obtain
samples from the posterior distribution, these samples could be used
to directly find the ``best,'' that is, most probable, motifs, obviating
the need for secondary analysis. We have obtained theoretical and
empirical results showing that the convergence of the Gibbs sampler is
even worse than previously realized. Our results reinforce the need
to convey the limitations of any estimates obtained using the Gibbs
sampler, and the need to develop more efficient Markov chain methods
for motif discovery.

Although our main result (Theorem~\ref{ThmSlowMix}) is phrased in
terms of a specific model, the methods used to prove this result are
very widely applicable to situations with i.i.d. data, where the data
are not necessarily generated according to the model, and where the
function $E \log f(X | \theta)$ is multimodal. The extent to which
slow mixing holds in other contexts will be determined by how
generally this multimodality condition holds, so we are currently
investigating this condition in detail.

\begin{appendix}\label{app}

%s7 #&#
\section{Bayesian asymptotics}

We quote a result from \citet{berk65} on Bayesian
asymptotics for i.i.d. observations. Let $f(x | \theta)$ be the
density (with respect to some $\sigma$-finite measure on a space
$\cY$) of each observation $X_i$ under the inference model,
parameterized by $\theta\in\Lambda$ where $\Lambda$ is a Borel
subset of a complete separable metric space. Let the true distribution
of the observations be
denoted by $G$. Define:\looseness=1
\begin{longlist}[(1)]
\item[(1)] The ``carrier'' of a distribution $P$: the smallest relatively
closed set having probability one under $P$.
\item[(2)] $P_n$: the posterior distribution of $\theta$ with $n$
observations, with
respect to a prior $P$ having carrier $\Lambda$.
\item[(3)] $\eta(\theta) \triangleq E\log f(X | \theta)$ where the
expectation is taken with respect to $ X \sim G$.
\item[(4)] $\eta^* \triangleq
\sup\{\eta(\theta)\dvtx \theta\in\Lambda\}$.

\item[(5)] $U_\delta\triangleq\{\theta\in\Lambda: \eta(\theta)
\geq\eta^*
-\delta\}$ for $\delta\geq0$.
\end{longlist}
Assume that:

\begin{longlist}[(1)]

\item[(1)] $f(x | \theta)$ is measurable jointly in $x$ and $\theta
$; for
$G$-almost every $x$, $f(x |\theta)$ is continuous in $\theta$.

\item[(2)] For all $\theta\in\Lambda$, $G\{x\dvtx f(x | \theta) > 0\}
= 1$.

\item[(3)] For any compact $F \subset\Lambda$, $E \sup_{\theta\in
F} |\log
f(X | \theta) | < \infty$.

\item[(4)] $\eta(\theta)$ is continuous.

\item[(5)] For any real number $r$ there is a co-compact set $D
\subset
\Lambda$ ($D^c = \Lambda\setminus D$ is compact) and a
cover $D_1,\ldots,D_K$ of $D$ such that
%
%e7.1 #&#
\begin{equation}
E \sup_{\theta\in D_k} \log f(X | \theta) \leq r, \qquad k \in\{1,\ldots ,K\}.
\label{EqnRegCond}
\end{equation}

\end{longlist}

With these assumptions, we have Theorem~\ref{ThmBerk}.
%
%th7.1 #&#
\begin{theorem}[{[\citet{berk65}]}] \label{ThmBerk}
For $G$-almost every sequence of observations $\{x_i\dvtx i \in\N\}$ and
any $\delta> 0$,
\[
\limsup_{n \rightarrow\infty} \biggl( \frac
{P_n(\Lambda\setminus U_\delta)}{P_n(U_\delta)} \biggr)^{1/n} \leq
e^{-\delta}.
\]
\end{theorem}
Theorem~\ref{ThmBerk} is a sub-result given in the proof of Berk's
main theorem. It says that the posterior probability of $U_\delta^c$
decreases exponentially in $n$. Here we have stated the result
slightly more generally than \citet{berk65} in the sense that we
replace his assumption (iii) with the only two relevant consequences of
that assumption: our assumptions (3)--(4). Also, in our assumption (5) we have
allowed a cover of $D$, whereas Berk's assumption (iv) takes $K=1$.
The extension to the case of general $K$ is immediate from his proof.

%s8 #&#
\section{Tools for bounding spectral gaps}\label{SecTools}

Let $P$ and $Q$ be transition kernels that are reversible with respect
to distributions $\mu_P$ and $\mu_Q$ on a (general) state space $\cX$
with countably-generated $\sigma$-algebra. Let $P|_B$ for $B \subset
\cX$ indicate the restriction of $P$ to $B$, which is defined to have
state space $B$ and transition probabilities identical to $P$ except
that any move to $B^c$ is rejected,
%
%e8.1 #&#
\begin{equation}
P|_B(x,D) = P(x,D) + \indic_{\{x \in D\}} P \bigl(x,
B^c \bigr), \qquad x \in B, D \subset B.\label{EqnTransRest}
\end{equation}
Also let $\mu_P|_B$ be the
restriction of $\mu_P$ to $B$, that is,
%
%e8.2 #&#
\begin{equation}
\mu_P|_B(dx) \triangleq\mu_P(dx) /
\mu_P(B), \qquad x \in B. \label{EqnDistRest}
\end{equation}
Then $P|_B$ is reversible w.r.t. $\mu_P|_B$.

For a partition $\{B_j\}_{j=1}^J$ of $\cX$, let $\bar{P}$ be the
\textit{projection matrix} of $P$ with respect to $\{B_j\}_{j=1}^J$,
defined to have state space $\{1,\ldots,J\}$ and\vspace*{1pt} $i,j$ element
equal to the probability that $P$ transitions to $B_j$, given that the
current state is in $B_i$. That is,
\[
\bar{P}(i,j) \triangleq\int\mu_P|_{B_i}(dx)
P(x,B_j),\qquad  i,j \in\{1,\ldots,J\}.
\]
The matrix $\bar{P}$ is reversible w.r.t. $\bar{\mu}$, where $\bar
{\mu}(j) \triangleq\mu_P(B_j)$.

%le8.1 #&#
\begin{lemma}[{[Madras and Zheng (\citeyear{madrzhen2003})]}]\label{ThmPowers}
For any $N \in\N$ we have\break
$\Gap(P) \geq\frac{1}{N} \Gap(P^N)$.\vadjust{\goodbreak}
\end{lemma}
Although \citet{madrzhen2003} state this result for finite state
spaces, their proof also holds for general state spaces.

%th8.1 #&#
\begin{theorem}[{[Madras and Randall (\citeyear{madrrand2002})]}]\label
{ThmDecomp}
$\!\!\!$Let $\mu_P = \mu_Q$, and let $\{B_j\}_{j=1}^J$ be any partition of
$\cX$. Assume that $P$ is nonnegative definite and let
$P^{1/2}$ be its nonnegative square root. Then
\begin{eqnarray*}
\Gap \bigl(P^{1/2} Q P^{1/2} \bigr) &\geq&\Gap(\bar{P})
\min_j \Gap( Q|_{B_j}),
\\
\Gap(P) &\leq&\Gap(\bar{P}),
\end{eqnarray*}
where $\bar{P}$ is the \textit{projection matrix} of $P$ with respect to
$\{B_j\}_{j=1}^J$.
\end{theorem}
%
%th8.2 #&#
\begin{theorem}[{[E.g., Sinclair (\citeyear{sinc1992})]}]\label{ThmConduct}
For $\cX$ finite define
\[
\Phi_{P} \triangleq\min_{B \subset\cX: 0<\mu_P(B)<1}\Phi_{P}(B)
, \qquad \Phi_{P}(B) \triangleq\frac{\sum_{x \in B} \mu_P(x)
P(x,B^c)}{\mu_P(B) \mu_P(B^c) }.
\]
Here $\Phi_{P}$ is
called the ``conductance,'' and $\Phi_{P}(B)$ is referred to as
the conductance of the set $B$. Then $\Gap(P) \leq2
\Phi_{P}$.
\end{theorem}

%th8.3 #&#
\begin{theorem}[{[Diaconis and Saloff-Coste (\citeyear
{diacsalo1996})]}]\label{ThmProdChain}
Take any $N \in\N$, and let $P_k$, $k = 0, \ldots, N$, be $\mu_k$-reversible transition
kernels on state spaces $\X_k$. Let $P$ be the
transition kernel with state $\bx= (x_0, \ldots, x_N)$ in the
space $\X= \prod_k \X_k$, given by
\[
P(\bx,d\by) = \sum_{k=0}^{N}
b_k P_k(x_{k}, dy_{k})
\delta_{\bx_{[-k]}}(\by_{[-k]})\, d\by_{[-k]}, \qquad \bx,\by\in\X
\]
for some set of $b_k > 0$ such that $\sum_k b_k = 1$, where $\delta$ is
Dirac's delta function, and where $\bx_{[-k]}$ indicates the vector
$\bx$
excluding $x_k$. $P$ is called a product chain with
``component'' chains $P_k$. It is
reversible with respect to $\mu_P(dx) = \prod_k \mu_k(dx_{k})$,
and
\[
\operatorname{Gap}(P) = \min_{k = 0, \ldots, N} b_k \operatorname{Gap}(P_k).
\]
\end{theorem}

Lemma~3.2 of \citet{diacsalo1996} states Theorem~\ref
{ThmProdChain} for finite
state spaces; however, the proof holds in the general case.
%
%le8.2 #&#
\begin{lemma} Take finite $\cX$ and
$\mu_P = \mu_Q$. If $\exists b > 0$ such that $bQ(x, y) \leq P(x,
y)$ for every $x,y \in\X$ such that $x \neq y$, then
$b\operatorname{Gap}(Q) \leq\operatorname{Gap}(P)$.
\label{moreTransitionLemma}
\end{lemma}
\begin{pf} The proof is nearly identical to that of Lemma~5.1 in \citet
{woodschmhube09a}.
\end{pf}

Lemma~\ref{moreTransitionLemma} is closely related to
Peskun ordering results; cf. \citet{pesk73,tier98,mira01}.

%th8.4 #&#
\begin{theorem}[{[\citet{sinc1992}, \citet{diacstro1991}]}]\label{ThmPathBound}
For $\cX$ finite, define a simple path
$\gamma_{x,y}$ between every ordered pair $x,y \in\cX$ in the graph of
the Markov chain with transition matrix $P$. A simple path is a
sequence of connected edges with no repeated vertices. Define the quantity
\[
\rho\triangleq\max_{(z,v) \in\E} \frac{1}{\mu_P(z)P(z,v)} \sum
_{\gamma_{x,y}\ni(z,v) } \mu_P(x)\mu_P(y) \operatorname{len}(
\gamma_{x,y}),
\]
where $\E$ is the set of edges,
where $\gamma_{x,y} \ni(z,v)$ is a path using the edge
$(z,v)$, and where $\operatorname{len}(\gamma_{x,y})$ is the number of edges in
$\gamma_{x,y}$. Then $\Gap(P) \geq\rho^{-1}$.
\end{theorem}
\end{appendix}

\section*{Acknowledgments}
The authors would like to thank Krzysztof Latuszynski
for assistance with one of the proofs, and the referee and
Associate Editor for their excellent feedback.

\begin{supplement}[id=suppA]
\stitle{Supplemental article}
\slink[doi]{10.1214/12-AOS1075SUPP} %[doi,text={...}] - jei reikia
%suskaldyti doi
\sdatatype{.pdf}
\sfilename{aos1075\_supp.pdf}
\sdescription{Provides additional proofs.}
\end{supplement}

% imsref loaded by akundreckaite, 2013-01-30 12:24:06
% imsref loaded by akundreckaite, 2013-01-30 12:24:32
%

\printaddresses

\end{document}